%% file: main.tex
\documentclass[graybox]{svmult}

\usepackage{rotating}
\usepackage{type1cm}        
\usepackage{comment}
  \usepackage[lofdepth,lotdepth]{subfig}
  
\usepackage{amsmath,amsthm}
\usepackage{makeidx}         
\usepackage{graphicx}        
\usepackage{multicol}        
\usepackage[bottom]{footmisc}

\usepackage[a4paper, total={6in, 8in}]{geometry} 
\usepackage[skip=10pt plus1pt, indent=40pt]{parskip} 
\usepackage{mdwlist}
\usepackage[printonlyused,withpage]{acronym}
\usepackage{makecell}
\usepackage{newtxtext}       %
\usepackage[varvw]{newtxmath}       
\usepackage{float}
\newcommand{\R}{\mathbb R}

\newcommand{\bD}{\mathbf{D}}

\newcommand{\bd}{\mathbf{d}}

\newcommand{\cA}{\mathcal{A}}
\newcommand{\eref}[1]{(\ref{#1})}

\newcommand{\dcp}{\mathbf{D}_{cP}}
\usepackage{hyperref}
\hypersetup{
    colorlinks=true,
    linkcolor=blue,
    filecolor=magenta,      
    urlcolor=cyan,
    pdftitle={Overleaf Example},
    pdfpagemode=FullScreen,
    }

\def \A{\mathcal{A}}

\def \C{\mathcal{C}}
\def \vol{\mathrm{vol}}
\makeindex             

\begin{document}

\title*{Studying Morphological Variation: Exploring the Shape Space in Evolutionary Anthropology}
\author{Shira Faigenbaum-Golovin and Ingrid Daubechies}

\institute{Shira Faigenbaum Golovin \at Duke University, 140 Science Dr, Durham, NC, \email{alexandra.golovin@duke.edu}
\and Ingrid Daubechies \at Duke University, 140 Science Dr, Durham, NC \email{ingrid@math.duke.edu }}

\authorrunning{S. Faigenbaum-Golovin, I. Daubechies}
\titlerunning{Learning the Underlying Structure in a Shape Space}

%
\maketitle

\vspace{-20mm}
\abstract{
We present results of a long-term team collaboration of mathematicians and biologists. We focus on building a mathematical framework for the shape space constituted by a collection of homologous bones or teeth from many species. The biological application is to quantitative morphological understanding of the evolutionary history of primates in particular, and mammals more generally. Similar to the practice of biologists, we leverage the power of the whole collection for results that are more robust than can be obtained by only pairwise comparisons, using tools from differential geometry and machine learning. 
This paper concentrates on the mathematical framework. We review methods for comparing anatomical surfaces, discuss the problem of registration and alignment, and address the computation of different distances. Next, we cover broader questions related to cross-dataset landmark selection, shape segmentation, and shape classification analysis. This paper summarizes the work of many team members other than the authors; in this 
paper that unites (for the first time) all their results in one joint context, space restrictions prevent a full description of the mathematical details, which are thoroughly covered in the original articles.
Although our application is to the study of anatomical surfaces, we believe our approach has much wider applicability.}

\section{Introduction}
\label{sec:1}

Machine learning methods provide tools to tackle a wide array of questions raised by large and complex data sets, with a core focus on examining and using the similarities and the differences among samples within a dataset. In some applications, observed scattered high-dimensional data 'points' (each of which can be quite complex and have significant internal structure -- e.g. images, curves, surfaces) are sufficiently similar to their neighbors that
it is reasonable to view
them as lying on an unknown manifold that one wants to understand better. 

\begin{figure}[h]
	\centering
	\captionsetup[subfloat]{farskip=0pt,captionskip=0pt, aboveskip=0pt}
	\subfloat[][]{ \includegraphics[height =0.4\textwidth]{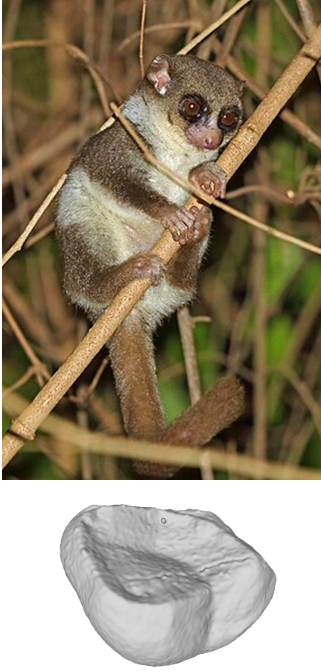} }  \hspace{0.8em}
	\subfloat[][]{ \includegraphics[height=0.4\textwidth]{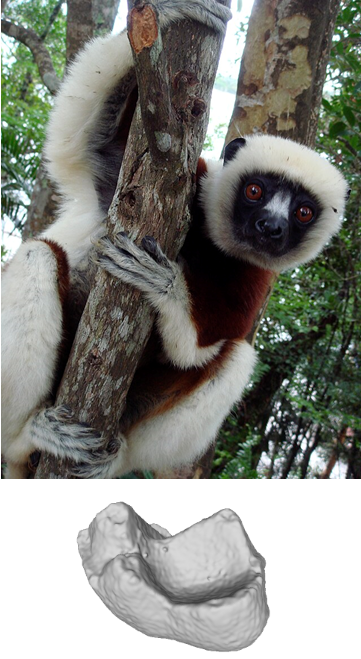} }\hspace{0.8em}
	\subfloat[][]{\includegraphics[height=0.4\textwidth]{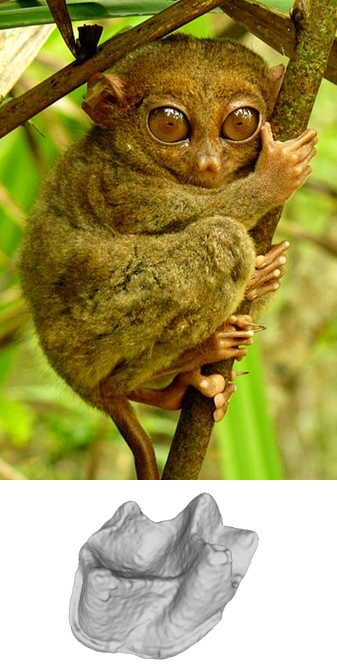}} \hspace{0.8em}
 \subfloat[][]{ \includegraphics[height=0.4\textwidth]{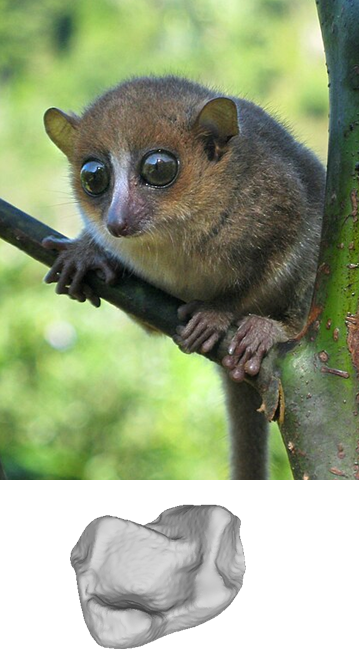} } 
	\caption{Comparison of dental morphology based on the diet. (A) Fruit eater; (B) Leaf eater; (C) Predator; (D) Omnivore (source: images wikimedia commons, meshes taken from morphosource  (AMNH:Mammals:M-80072, AMNH:Mammals:M-100832, AMNH:Mammals:M-196480, AMNH:Mammals:M-174537) \cite{morphosource}.
}
	\label{fig:teethAndDiet}
\end{figure}

In our particular application, we typically consider a collection of anatomical shapes, usually with some rough correspondence between the shapes, and we aim to study the underlying shape space. 
In the current paper, we survey how applying tools from the field of differential geometry and machine learning aids in studying anatomical surfaces, in the framework of Evolutionary Anthropology. Comparing such surfaces is vital for the understanding of the evolutionary history of primates; quantitative comparisons help in answering questions of interest in evolutionary anthropology 
Morphology and relative sizes of the anterior dentition have long played a role in the interpretation of systematics and
evolution of early primates \cite{simons1961dentition, rose1984gradual, rasmussen1995dentition}. Understanding the relationship between skeletal morphology, function, and ecology in this group is critical for understanding the evolutionary origins of the order to which humans belong. It is generally thought that particular tooth forms are suited to processing essential food items of respectively particular forms and material properties, and that mammalian teeth/dentitions are often unique taxonomic identifiers
for more details see~\cite{boyer2008relief}.  Examples of the different teeth of a fruit eater, a leaf eater, a predator (of arthropods), and an omnivore are shown in Figure \ref{fig:teethAndDiet}.
All the teeth used in this study can be downloaded from MorphoSource \cite{morphosource}.
Although the surfaces illustrated in Figure \ref{fig:teethAndDiet} differ considerably, tooth surfaces in a much wider collection acquired from many species show more gradual transitions, in particular between species closely related to each other; the variations are sufficiently incremental that the surfaces can reasonably be viewed as samples lying on a manifold.

In many cases, only a surface cast of the original tooth or bone is available in our collection; for this reason, we work with surfaces rather than volumetric data. The surfaces are scanned as scattered data and later represented as triangular meshes. 
Compiling the geometric information in these triangular meshes is a challenge, due to the nuanced variations in anatomical shapes. Our team follows a workflow that consists of several steps. First, pre-processing steps are performed, starting from denoising and cleaning the triangulated arrays
followed by a coarse alignment of these surfaces \cite{puente2013distances} as well as a first estimation of a robust and biologically meaningful distance for each pair of surfaces. Next, the variation between the shapes, using the power of {\em collection} of surfaces is addressed. A (up-to
-date only until the time of its writing) survey of several automated geometric morphometric methods (some also discussed below) can be found in \cite{gao2018development}.

Although surface analysis, such as comparison and registration, is extensively addressed in the computer graphics or Computer Aided Geometric Design (CAGD) communities, the algorithms developed there cannot simply be ported to our application. On the one hand, the features of interest in anatomical shapes and their variations are typically too subtle to be considered useful features for computer graphics methods; on the other hand, it turns out that leveraging manifold analysis to study a whole collection of surfaces leads to much better results. In addition, the research questions stemming from evolutionary biology introduce a unique perspective. 

\begin{figure}[h]
	\centering
	\captionsetup[subfloat]{farskip=0pt,captionskip=0pt, aboveskip=0pt}
	\subfloat[][]{ \includegraphics[width=0.08\textwidth]{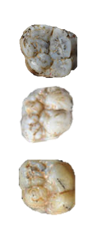} } \hspace{0.3em}
	\subfloat[][]{ \includegraphics[width=0.24\textwidth]{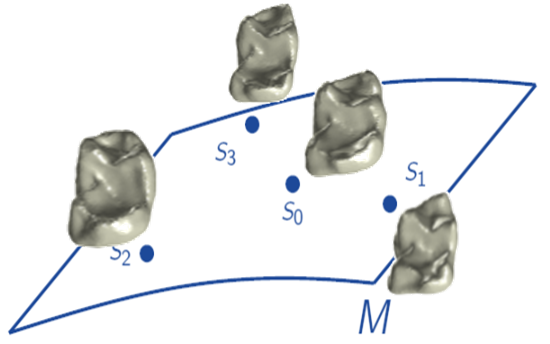} }\hspace{0.3em}
	\subfloat[][]{\includegraphics[width=0.19\textwidth]{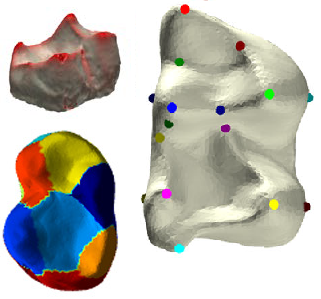}} \hspace{0.3em}
 \subfloat[][]{ \includegraphics[width=0.40\textwidth]{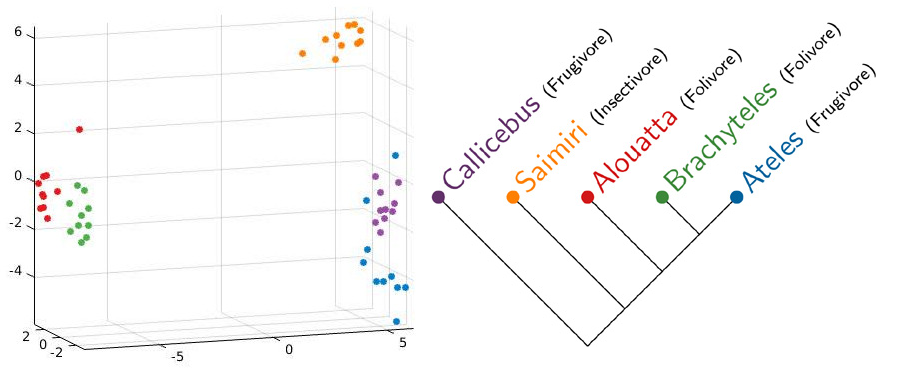} } 
	\caption{Working flow in a shape space. (A) Data acquisition; (B) Alignment, registration, and distance calculation on the manifold of shapes; (C) Automatic Shape analysis (e.g., from top left clockwise, curvature estimation (ariaDNE \cite{shan2019ariadne}), landmark sampling, consistent segmentation of shapes in the collection ; (D) Phylogeny of different species groups. Images adopted from \cite{shan2019ariadne, gao2015hypoelliptic}.}
	\label{fig:shapeSpaceAnalysisFlow}
\end{figure}

In this paper, we survey the research on anatomical surfaces carried out by a team of biologists and mathematicians over more than 15 years, involving a total of more than 20 researchers (at different times) focusing on the mathematical aspects of the work. We cover the methods that were developed for the comparison of the surfaces, discuss the problem of registration and alignment, and review the calculation of the different distances. Next, we cover broader questions related to cross-dataset landmark selection, shape segmentation, and shape classification analysis. We conclude this paper with our view of the future of the shape space realm. These steps are illustrated in Figure \ref{fig:shapeSpaceAnalysisFlow}. Our main analysis focuses on disc-type surfaces, such as the part of a tooth above the gum); most of our methods can be (and some have been) applied to sphere-type surfaces (which have no boundary). So far, we have not extended the work to higher genus surfaces (which have holes).

The remainder of the paper is organized as follows.
In Section \ref{sec:Data Acquisition} we review the data acquisition process, from CT scan to clean surfaces. Next, in Section \ref{sec:reg_landmark} we review the methods that result in robust distance metrics between the shapes, starting with alignment, landmarking, and registration. Then in Subsection \ref{sec:Distance} we detail several methods for distance calculation (among them is the Conformal Wasserstein local distance, continuous Procrustes distance, as well as the Horizontal-based diffusion distance). Quality control aspects are provided with the distance accuracy validation in Section \ref{sec:Quantifying Errors}. Next, we dive into the manifold learning analysis in Section \ref{sec:Manifold Learning}, where we briefly survey the classical tool for studying manifold geometry, Diffusion Map in Subsection \ref{sec:DM_definition}. Next, the theoretical and practical aspects of decomposing manifolds by exploring fibers and the base manifold \label{sec:manifold learning_Fiber} are explained. Later, we provide details about the preeminent approach of horizontal random walks and diffusion processes on fibre bundles \label{HDM_definition}. In Section \ref{sec:ariaDNE} we are finally ready to shed light on morphology via statistical analysis and machine learning methods. In Table \ref{tab:code} links to all the code discussed in this paper appear, while the available datasets used for the study are listed in Table \ref{tab:data}.


\textbf{List of abbreviations used in this survey}
\begin{itemize*}

\setlength{\itemsep}{0pt}
  \setlength{\parskip}{0pt}

    \item [] \textbf{MSE} - Mean Square Error
    \item [] \textbf{MST} - Minimum Spanning Tree
    \item [] \textbf{PCA} - Principal Component Analysis
    \item [] \textbf{Auto3DGM} - Automated 3D Geometric Morphometrics 
    \item [] \textbf{MSPE} - Maximum Mean Squared Prediction Error

    \item [] \textbf{DM} - Diffusion Map 
    \item [] \textbf{SGE} - Semi-group Error
    \item [] \textbf{cP} - Continuous Procrustes
    \item [] \textbf{cPDist/CPD} - Continuous Procrustes distance
    \item [] \textbf{cWn} - Conformal Wasserstein Neighborhood Dissimilarity Distance 
    
    \item [] \textbf{HDM} - Horizontal Diffusion Maps
    \item [] \textbf{HBD} - Hypoelliptic Diffusion Distance
    \item [] \textbf{HBDM} - Hypoelliptic Base Diffusion Maps
    \item [] \textbf{HBDD} - Hypoelliptic Base Diffusion Distance

    \item [] \textbf{DNE} - Dirichlet Normal Energy
    \item [] \textbf{ariaDNE} - a Robustly Implemented Algorithm for Dirichlet Energy of the Normal
    \item [] \textbf{EOP} - Eyes on the Prize
    \item []
\end{itemize*}

\vspace{-10mm}

\section{Data Acquisition and Pre-Processing}
 \label{sec:Data Acquisition}
\begin{figure}[h]
	\centering
	\label{fig:c}\includegraphics[width=0.8\textwidth,keepaspectratio]{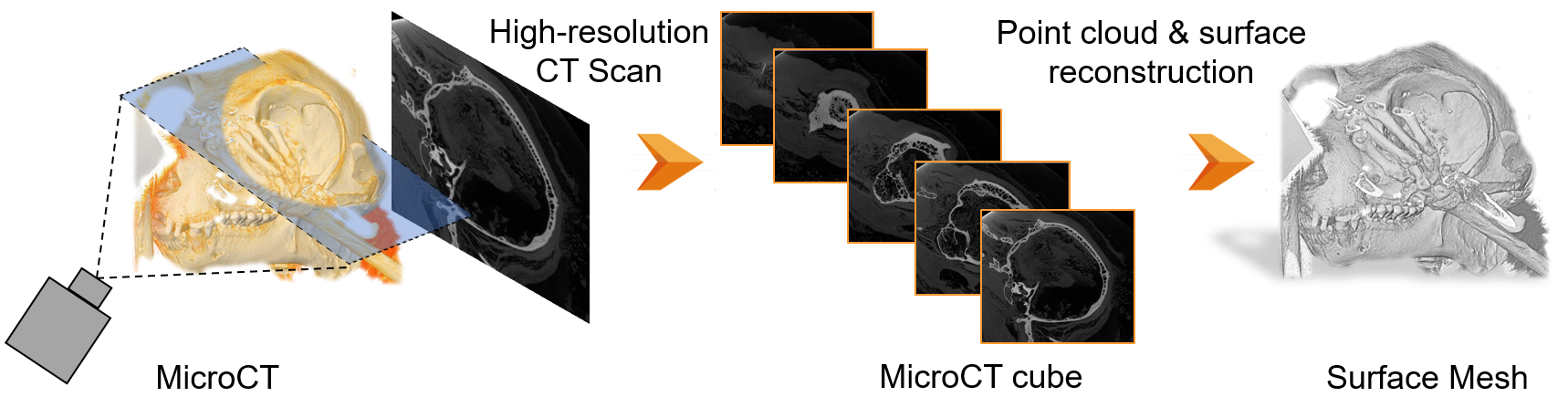} 
	\caption{Illustration of the scanning process of a skull: from CT scan to MicroCT cube, to point cloud and surfaces reconstructed and cleaned.\protect\footnotemark}
	\label{fig:acquisition and denoising}
\end{figure}

We shall visit two types of example surfaces: disk-type surfaces obtained by scanning with high resolution the portion above the gumline of second mandibular molars of prosimian primates and nonprimate close relatives, and sphere-type surfaces obtained by scanning talus bones of primates and close relatives. 
Digitized surfaces were obtained from high-resolution X-ray computed
tomography ($\mu$CT) scans, 
using the 
A Nikon XT H 225 ST $\mu$CT machine at Duke University’s Shared Materials Instrumentation
Facility (SMIF). 
For greater detail on the protocol followed, see e.g. the description in \cite{yapuncich2019digital}. Figure \ref{fig:acquisition and denoising}  illustrates the steps to acquire a surface mesh by scanning a (small) real object using an X-ray $\mu$CT system.

The MicroCT cubes were used to create a 3D digital point cloud collection using 
Avizo or dragonfly software \cite{Avizo, dragonfly}, after which a mesh was created using via 
standard triangulation algorithms. 
The procedure of converting the CT scan images to point clouds and meshes is typically prone to many errors, leading to  inconsistencies in the meshes;
among such artifacts are
incorrect identification of nearest neighbors of certain points, reversed neighboring connectivity, resulting in failure to generate a smooth surface (with some faces flipped), spikes, and holes. These are then manually identified, and off-the-shelf software GeoMagic  
is used to address these issues. 
All scans are publicly available on MorphoSource \cite{boyer2016morphosource}, an online repository specifically designed to archive 3D data.
 
\footnotetext{CT scan images courtesy of The Duke Lemur Center that provided access to these data, originally appearing in Yapuncich et al. (2019), the collection of which was funded by NSF BCS 1540421 to Gabriel S. Yapuncich and Doug M. Boyer. The files were downloaded from www.MorphoSource.org, Duke University.}

\vspace{-10mm}
\section{From Scattered Surfaces to Metric Shape Space}

\label{sec:reg_landmark}
Once the shapes are scanned, the point clouds extracted, and surfaces are created and cleaned, we are ready to dive into the shape analysis and the comparison question. The main challenges for calculating the similarity between shapes are alignment, registration, and distance calculation. 
Alignment and registration are tightly connected but represent different tasks. While \textit{alignment} moves one object into alignment with another object so that they both share the same coordinate system, \textit{registration} is the mapping from one surface to the other that connects areas with corresponding function and geometry.

Calculating the distance between two registered objects typically amounts to integrating (or summing), over one of the objects, the distance of each point to its corresponding point on the other object. On the other hand, the optimal correspondence map from one object to the other is usually found by optimization, 
minimizing the value of the distance between the objects. Pairwise registration and distance computation thus go hand-in-hand. In this section, we review several algorithms used by our team for these tasks. Often our analysis workflow uses several different approaches, with simple algorithms at first, and more sophisticated ones at a later stage.

\subsection{Alignment, Landmarking, and Registration}
\label{sub_sec_regi_align}

\begin{figure}[h]

 \centering
	\captionsetup[subfloat]{farskip=0pt,captionskip=0pt, aboveskip=0pt}
	\subfloat[][]{ \includegraphics[width=0.4\textwidth]{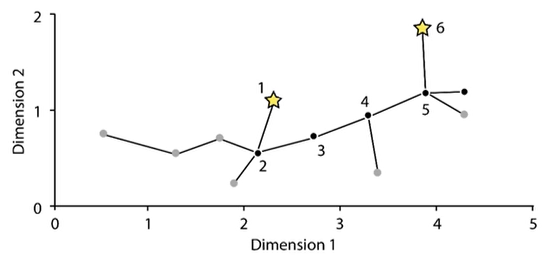} } 
	\subfloat[][]{ \includegraphics[width=0.6\textwidth]{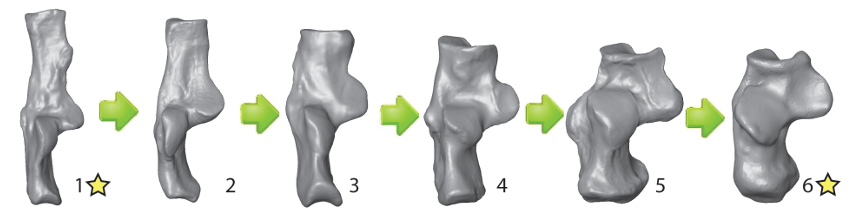} }
	
	\caption{Alignment of two different shapes labeled as 1 and 6, by composing a series of incremental different bones (right) found via the MST (left) \cite{boyer2015new}.}.
	\label{fig:alignment}
\end{figure}

\textbf{Alignment:} The goal of shape alignment is to establish a meaningful correspondence between corresponding points, features, or structures on different shapes, ensuring that they are spatially aligned or registered with each other. For computational purposes, the meshes are first down-sampled prior to analysis. It is useful to note  that finding a joint coordinate system is easier for similar shapes than for different shapes. This idea is incorporated into the Automated 3D Geometric Morphometrics (Auto3DGM) tool \cite{Auto3dgm}, one of the earliest tools developed by our team. It first uses PCA to compute the three principal axes for each shape. For each pair of shapes, an alignment is carried out based on best matching their principal axis, after which the associated Procrustes distance between the shapes is calculated. (This distance is defined later in \eqref{eq:cont_P_dist}). These are only  intermediate steps -- the alignment and distance obtained in this way are typically not very good when the two shapes are too dissimilar. However, the
approximate distances can be used to determine a minimum spanning tree (MST) for the collection; between two elements that are close to each other, the alignment is typically quite good. By aligning shapes with each other sequentially, following the MST, even very dissimilar bones can be aligned correctly, using intermediate bones in the process, as illustrated in Figure \ref{fig:alignment}.  More details can be found in \cite{boyer2015new}). 

\textbf{Registration and Landmarking:} Registration is the mapping of every piece of one shape to another while taking into account the corresponding features. Shape matching is an area of active research in computer graphics \cite{deng2022survey, van2011survey} and also the prior art appearing in \cite{ravier2018algorithms}. In order for the registration to be useful for our
biological collaborators, it must 
(a) be consistent with the field knowledge about the anatomically equivalent (biologically homologous parts); (b) take into account the variation present in other species of the sample; and (c) correspond to a {\em smooth} map between shapes: 
the ability to compute high-quality diffeomorphisms between shapes would allow the use of a wider range of statistical techniques in order to compare
shapes in a collection.

One plausible solution is to choose landmarks that capture some geometric features of the surfaces. 'Landmarking by experts' is a manual process (now typically carried out on a computer screen) where an equal number of geometrically or semantically meaningful feature points, the \emph{landmarks}, are marked on a surface. The selected landmarks are certified by domain experts to be in consistent one-to-one correspondences. The methodology of landmark selection emphasizes a comprehensive and balanced decision between sharp geometric features and points taken farthest away from other points under certain metrics (for a discussion about the popular methods see \cite{watanabe2018many, warmlander2019landmark}). While experts can manually suggest the important feature points, and put into correspondence areas of high and low curvature, they have a harder time prescribing the correspondence map in intermediate regions, between these feature points. In rare cases, experts can also disagree on landmark correspondences. The imposition that the same exact number of landmarks be picked on all the surfaces in the collection can also be a limiting aspect. Therefore, we considered several automatic procedures for locating the important points, described below. In our workflow at present, we typically start with a simple pairwise approach the results of which provide the point of departure for more sophisticated methods, used further
in the workflow.

\begin{figure}[h]

\centering
	\captionsetup[subfloat]{farskip=0pt,captionskip=0pt, aboveskip=0pt}
	\subfloat[][]{\includegraphics[width=0.48\textwidth]{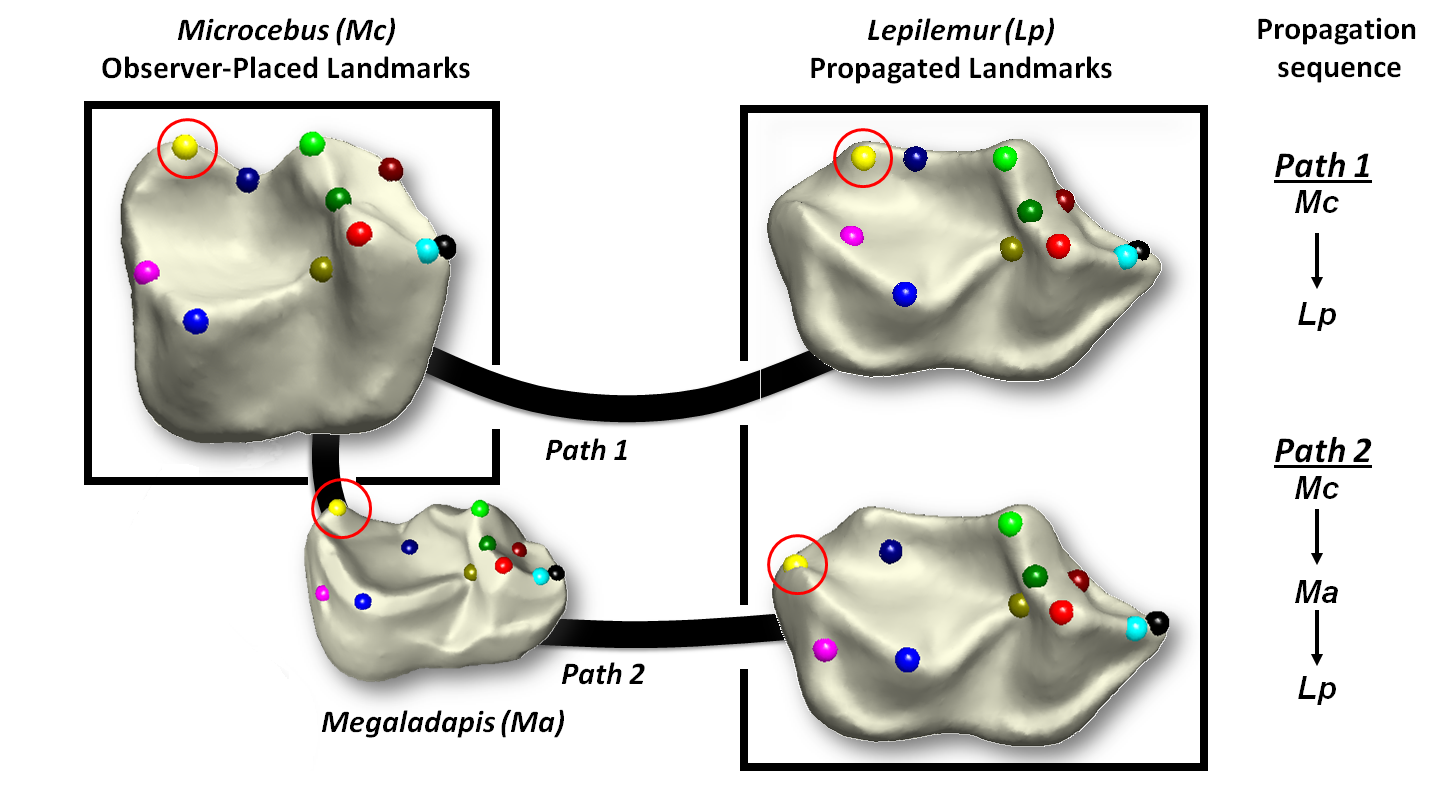} } \hspace{1em}
 \subfloat[][]{\includegraphics[width=0.48\textwidth]{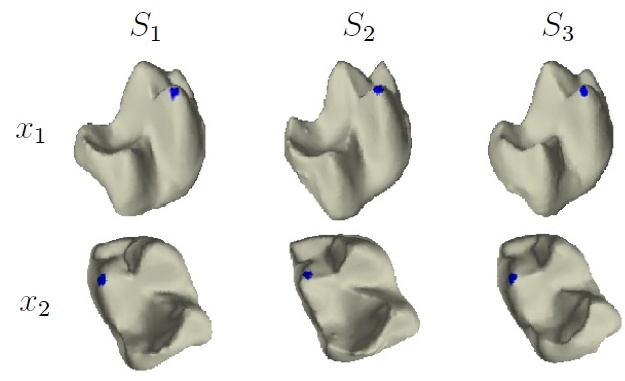} }
 
	\caption{(A) Landmark propagation via continuous Procrustes distance. Propagating from Microcebus to Lepilemur is not equivalent if their ancestor Megaladopis is also used  \cite{boyer2011algorithms}; (B) Landmarking via Matching Pursuit algorithm. Points that are in correspondence with each other using an eigenvector (each row is a correspondence to a different eigenvector) \cite{gao2015hypoelliptic}.}.
	\label{fig:landmarking}
\end{figure}

\textbf{Continuous Procrustes Landmarking:} In the discrete implementation of the continuous Procrustes distance (defined below in equation \eqref{eq:cont_P_dist}), one starts by distributing equal numbers of approximately evenly spaced points on two surfaces; an iterative procedure then determines the correspondence between the two surfaces that minimizes the Procrustes distance between the two point sets. This provides us with an out-of-the-box landmarking. However, since this landmarking is performed pairwise, the propagation of landmarks between shapes is not always consistent, and can change if a morphologically close shape is added: in the left panel of Figure \ref{fig:landmarking} landmark propagation from shape A to B changes if we first propagate from A to C and then from C to B.

\vspace{-8mm}
\textbf{Matching Pursuit Algorithm Landmarking:}
An alternative method for landmark selection was proposed in \cite{gao2015hypoelliptic}, using a
Matching Pursuit algorithm. It exploits the similarity (or rather dissimilarity) between surfaces, and among different portions of the surfaces, to create a large matrix, $L_{\alpha, *}^H$, the graph hypoelliptic Laplacian, generalising the graph Laplacian, on the union of all the sample points of all the surfaces. An exact correspondence map between the surfaces (normalized to all have area 1) would lead to a sparse eigenvector of that matrix. This suggested the following optimization approach to find a vector $x$ that would (approximately) encode corresponding landmarks:
$$ min_{ ||x||_2 = 1} \,
\left[ \,||L_{\alpha, *}^H x||_2^2 + \lambda ||x||_1 \, \right]$$

More details can be found in in \cite{gao2015hypoelliptic}. 
The algorithm produces corresponding points on the entire collection of surfaces, even though this is not built into the functional {\em a priori}. This method was the first, developed by our team, that leverages the power of considering a {\em collection} of surfaces. See Figure \ref{fig:landmarking} (B) for an illustration.

\begin{figure}[h]

 \centering
	\captionsetup[subfloat]{farskip=0pt,captionskip=0pt, aboveskip=0pt}
	\subfloat[][]{ \includegraphics[width=0.15\textwidth]{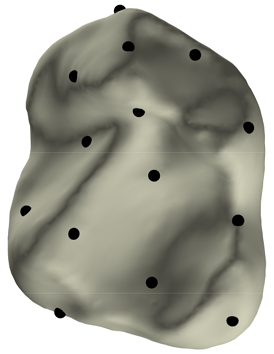} }  \hspace{0.5cm}
	\subfloat[][]{ \includegraphics[width=0.17\textwidth]{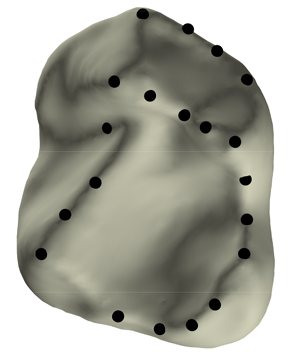} } \hspace{0.5cm}
 \subfloat[][]{ \includegraphics[width=0.35\textwidth]{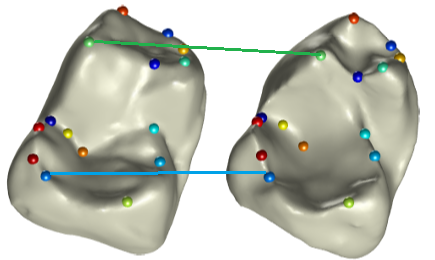} }
	
	\caption{Shape registration via automatic Landmarking \cite{gao2019gaussian}. (A)-(B) demonstrate different landmarking methods: geodesic farthest point sampling in (A), by Gaussian process landmarking in (B) . (C) After Gaussian process landmarking on each of the two surfaces, bounded distortion matching identifies corresponding pairs of landmarks (spheres with matching colors); note that each surface may have many landmarks (not shown) that
    do not correspond to a corresponding landmark on the other surface. }.
	\label{fig:landmarkSelection_Gao}
\end{figure}

\textbf{Gaussian Process Landmarking:}
A different approach to the 'automatic' identification of salient
points, which can be interpreted as landmarks, was introduced in \cite{gao2019gaussian, gao2019gaussian_2}. It is essentially a 
Kriging method, based on successively selecting points of maximum mean squared prediction error (MSPE) with respect to the Gaussian process on the anatomical surface, with a variance-covariance structure determined by the heat kernel of the surface, viewed as a Riemannian manifold; the landmarks are then selected successively, each time picking a new landmark as the point on the surface with the largest variance conditioned on all the previously selected landmarks. Figure \ref{fig:landmarkSelection_Gao} compares two methods for landmark selection: produced by geodesic farthest point sampling (A) \cite{moenning2003fast}, and by Gaussian processing landmarking (B) \cite{gao2019gaussian}. As can be seen, The Gaussian sampling (B) prioritizes choosing semantically meaningful features for geometric morphometrics and agrees better with points selected by biological experts.

After landmarks are identified on two surfaces $S_1$ and $S_2$, the best registration via landmark matching is determined. Among all the candidates obtained by Gaussian process landmarking on each of the surfaces, 
a maximal subset of geometrically consistent correspondences is determined; next interpolation is used found to obtain a diffeomorphism between the surfaces $S_1$ and $S_2$ that connects the landmarks in this subset; this results in two fully registered surfaces. See Figure \ref{fig:landmarkSelection_Gao} (C) for an illustration of landmark matching. MATLAB code for the surface registration algorithm is publicly available at \cite{GPLmkBDMatch}.

\vspace{-10mm}


\textbf{Consistent registration methods:}
As pointed out above, in the discussion of Auto3DGM Alignment, propagating landmarks via best correspondences that are determined by pairwise comparisons of surfaces (whether done
via the 
Continuous Procrustes approach, identifying best possible rigid rotations and pairwise mapping, or via e.g. Gaussian landmarking), following paths on the corresponding distance graph for a whole collection by determining the Minimum 
Spanning Tree, can lead to inconsistencies (see 
Figure \ref{fig:landmarking}A, as well as the discussion below about the accuracy of the distance induced by the Continuous Procrustes process (Section \ref{sec:Quantifying Errors}). To avoid this sensitivity to the selection of surfaces included in the study, we have also introduced alternative methods that use 
diffusion networks to propagate information, or to define a consistent reparametrization of the full collection. We shall come back to these methods below, after discussing different distances defined on collections of surfaces.

\subsection{Distances in the Shape Space}
\label{sec:Distance}
We turn now to the question of comparing the shapes qualitatively, studying them statistically, and testing different hypotheses of evolution. Quantifying and comparing complex shapes is a key component of 
evolutionary morphology.
The standard approach used by biological morphologists computes the {\em Procrustes distance} between surfaces, based on the preliminary identification of homologous features (corresponding points or {\em landmarks} on each of the surfaces in the collection to be studied, based on the domain knowledge of the landmarking expert of the functionality of different geometric features of the bones or teeth). 
The collection of the landmark points on each surface can then serve as a proxy for its surface, and 
distances between shapes are computed based on 
the following definition.

\begin{definition} [Procrustes distance]
    The \emph{Procrustes distance} between two surfaces $S_1$, $S_2$ is computed using the following procedure:
\begin{enumerate}
\item\label{item:3} Specify two sets of (operationally) homologous landmarks $\left\{ x_{\ell}^{\left( 1 \right)}\mid 1\leq \ell\leq L \right\}$ and $\left\{ x_{\ell}^{\left( 2 \right)}\mid 1\leq \ell\leq L \right\}$ on $S_1$ and $S_2$, respectively.
\item\label{item:4} Compute the distance between $S_1$ and $S_2$ by minimizing the energy functional
  \begin{equation}
    \label{eq:cP-variational}
    d_{\mathrm{cP}}\left( S_1,S_2 \right)=\inf_{T\in\mathbb{E}\left( 3 \right)}\left(\frac{1}{L}\sum_{\ell=1}^L \left\| T\left(x_{\ell}^{\left( 1 \right)}\right)-x_{\ell}^{\left( 2 \right)} \right\|^2\right)^{\frac{1}{2}}
  \end{equation}
where $\mathbb{E} \left( 3 \right)$ is the group of rigid motions in $\mathbb{R}^3$.
\end{enumerate}
\end{definition}

 The main disadvantage of this method is that it requires  corresponding landmarks as input; pre-determining these 
 can be a time-consuming process. 
 
 One of the first tasks our team was asked to tackle, in our collaboration with biological morphologists, was to define 
 distances (and algorithms to compute them) that would not require preliminary expert-defined landmarking, but that would nevertheless be as useful for quantification of similarities and differences in their biological studies as their usual Procrustes distances. In our first approach, we used mass transportation techniques to define 
 the {\em Conformal Wasserstein neighborhood distance}. In a first step,  the disk-type surfaces are mapped conformally to a disk; the conformal factor for each is then interpreted as a "landscape" on the disk to which we apply a mass transportation principle. (A similar approach is also possible for sphere-type surfaces, mapping them conformally to a sphere instead of a disk.) To define a mass transportation 
 distance between the two landscapes, we need to introduce a 
 (dis)similarity metric distance between points $v$, $w$ in the disk adapted to our problem: 
 to make sure that the resulting distance
 between the surfaces is independent of the (non-unique) conformal flattenings that were performed first, we define this metric distance in a conformally invariant way. The distance between points $z$ (in the landscape determined by $\mu$) and $w$ (in the landscape for $\nu$) is set to be the $L^1-$distance between the landscapes around the points, restricted to a small local hyperbolic neighborhood of each, after moving them via a M\"obius transformation (i.e. a
 conformal map on the disk) to the center of the disk. 
 This pointwise distance is then used to define a Wasserstein distance \cite{rubner2000earth} between the two configurations. \cite{lipman2011conformal}. We have thus:

\begin{definition} [Conformal Wasserstein neighborhood distance (cWn)] Given two disk-type surfaces $S_1$, $S_2$
we map each conformally to the hyperbolic disk $D = \{z | |z| < 1\}$, with corresponding conformal factors $\mu$, $\nu$. For a fixed choice of $R<\!<1$, we define the \emph{conformal neighborhood dissimilarity distance} $\bd^R_{\mu,\nu}(v,w)$
as the minimum over $\theta \in [0,2\pi]$ of the 
$L^1-$distance (with respect to the hyperbolic surface area of the disk $D$) between the restrictions to the disk $N(0,R)=\{\,z\,;\,|z|<R \, \}$ of the functions $\mu \circ {\tau_1}^{-1}$ and $\nu \circ {\tau_2}^{-1}\circ R_\theta$, where $\tau_1$, $\tau_2$ are, respectively, M\"obius transformations that map $v$ and $w$ to the center 0 of $D$, and $R_\theta$ is the rotation around 0 by an angle $\theta$.  Equivalently, $\bd^R_{\mu,\nu}(v,w)$ is defined as:
$$
\bd^R _{\mu, \nu}(v,w)=\!\!\inf_{m \in \mathcal{M_D}, m(v)=w}\!
\left[\,\int_{N(v,R)}|\mu(z)-(m*\nu)(z)\,|\,d_{vol_H}(z)\right],
$$

where $\mathcal{M_D}$ is disk-preserving Mobius group, $N(v,R)= {\tau_1}^{-1}\left(N(0,R)\right)$, and $d_{vol_H}$ is the volume form for the hyperbolic disk.

The distance between the two surfaces $S_1$ and $S_2$ is then defined by the standard Kantorovich approach to mass transportation distances:

\begin{equation}
\label{eq:conf_Wass_diss_dist}
\mathcal{D}^R_{\mbox{\small{cWn}}}(S_1, S_2)\,=\,
\inf_{\pi \in \Pi(\mu,\nu)} 
\int_{\bD{\times}\bD}
\bd^R_{\mu, \nu}(v,w)\,d\pi(v,w)~,
\end{equation}

where $\Pi(\mu,\nu)$ is the set of all probability measures 
on ${\bD{\times}\bD}$ of which the marginals coincide with respectively $\mu$ and $\nu$. 

\end{definition}

An attractive feature of this distance is that it determines
a (dis)similarity between points on $S_1$ and $S_2$ by
comparing their respective surroundings on each surface, which is quite local; on the other hand, it also 
takes into account the global setting, by integrating over the whole surface.




In practice, the cWn distance, with its multiple steps each involving optimization (including one over the 5-dimensional 
set of M\"obius transformations), was fairly heavy computationally. We therefore also introduced another approach, the {\em Continuous Procrustes distance} (cP), \cite{al2013continuous}, inspired by the Procrustes distance but landmark-free. 

When landmarks are given on a collection of surfaces, one could
conceivably interpolate between them and thus use this to define
correspondence maps defined on the full surfaces (and not just on the landmark points), or, in a computational setting, on a much larger number of discrete points than only the landmarks. 
In the absence of landmark points between which to interpolate, a reasonable approach is to pick a large number $N$ of "reasonably well distributed" points on each surface and view bijections between the $N-$tuples on different surfaces as discrete approximations to possible correspondence maps. If "reasonably well distributed" means that the Voronoi region around a random sample point should have an area of (approximately) $1/N$ (assuming the total area of each surface has been normalized to 1), then these bijections are discrete stand-ins for area-preserving maps. This motivated the following 
definition of the {\em Continuous Procrustes distance}:


\begin{definition} [Continuous Procrustes distance between surfaces (cP)] Given two surfaces $S_1$, and $S_2$, we
denote by
$\cA(S_1, S_2)$  the set of all area-preserving \cite{al2013continuous} diffeomorphisms, and let $\mathbb{E}_3$ be the Euclidean group on $\R^3$. Then the Continuous Procrustes distance between $S_1$ and $S_2$ is given as

\begin{equation}
\dcp(S_1,S_2)= \inf_{\C\in \A(S_1,S_2)} \,\left[\,
\inf_{R \in \mathbb{E}_3}\,\left(\, \int_{S_1} \, ||R(x)-\C (x) ||^2 \,d\vol_{S_1}(x) \, \right)^{1/2}
\, \right]\,.
\label{eq:cont_P_dist}
\end{equation}
\end{definition}

Note that the cWn-distance was \emph{intrinsic}: it used 
only
information ``visible'' from within each surface, such as geodesic
distances between pairs of points; consequently it would not be able to distinguish
a surface from any of its isometric embeddings in 3-D. (For instance, it would not distinguish a Swiss roll 2D-surface from a rectangle.) The cP uses {\em extrinsic} information about the surfaces; its computation also provides, as a side result, a mapping between the surfaces, which is the continuous Procrustes
registration described in the previous section. 

Comparing the results of the cWn and cP distances with the Procrustes distances derived from expert-placed landmarks (which to our biological collaborators served as "ground truth") showed
that the cP distance performs better \cite{boyer2011algorithms}; 
we have therefore used it more systematically, at least in its
numerical implementation. Moreover, we proved in \cite{al2013continuous} that when the cP distance between two surfaces is small, the corresponding minimizing area-preserving map is also close to the minimizing conformal map; this has the important consequence that in a numerical implementation, the 
(practically infeasible) search over area-preserving maps can be replaced by a (much more manageable) search over conformal maps, with a guarantee on the error on the estimate for the cP distance made by substituting this easier optimization. 

The cP correspondence maps are determined pairwise, and turn 
out (unsurprisingly) to be non-transitive, which is a disadvantage for collections of bone or tooth shapes from species that are sufficiently closely related that biologists expect transitive correspondence maps between all of them. Uniqueness of correspondence could be restored by considering
only those maps associated with a Minimizing Spanning Tree in the distance graph, but this then suffers from the sensitivity 
to inserting extra specimens in the collection pointed out before. Finally, as reported in \cite{boyer2012earliest} there exist anomalies in some of the pairwise correspondence maps between specimens. 
Thus, it is not clear how the resulting cP correspondence maps could be incorporated in a more traditional geometric morphometric workflow. 

In a later subsection, we introduce a new mechanism that uses tools from differential geometry, that uses the cP correspondence maps (or other correspondence maps, if available) as input to a reparametrization that leads to higher-quality, transitive correspondence maps and an improved distance graph. In order to define this new distance, we first need to introduce some extra mathematical machinery. We will discuss Diffusion maps, the corresponding classical Diffusion Distances, and finally Horizontal-based diffusion in full length in  Section \ref{HDM_definition}. However, for the sake of completeness in our list of distances, we provide here  place-holding definitions for Diffusion distance and Horizontal-based diffusion distance and refer the reader to additional details in Sections \ref{sec:DM_definition} and \ref{HDM_definition}. 

\begin{definition} [Diffusion Distance (DD)]
Given a collection of objects, endowed with a distance metric.
Let $\left(W_t\right)_{t > 0}$ be the diffusion matrices constructed by means of these distances, as defined in \eqref{Diffusion_matrix}, and let 
$\psi_1$, $\ldots, \,\psi_L$ be the first few dominant non-trivial eigenvectors of $W_{\tau}$, for one appropriately chosen $\tau >0$ (see below), 
with corresponding eigenvalues $\lambda_1$, $\ldots,\, \lambda_L$. Then we define a parametrization of the objects by means of the values of these eigenvectors (i.e. object $k$ is parametrized by the values $\psi_{\ell}(k)$, with $\ell$ ranging over $1,\, \ldots,\, L$) and for each $\gamma >0$ we define a $\gamma$-dependent diffusion map $D_{\gamma}$ by setting 
\begin{equation}
  \label{eq:DD}
D_{\gamma}(k,k')^2=\sum_{\ell=1}^L \lambda_{\ell}^{\gamma} 
\,|\psi_\ell(k) - \psi_\ell(k')|^2\,.
\end{equation}
\end{definition}
The Diffusion Distance can be (and has been) used for much 
more general collections than the bones and teeth we consider here; in our case, we shall apply the Diffusion distance approach to approximative cP distances between bone or teeth 
surfaces that we compute from "well-distributed" samples on those surfaces. We shall see in Section \ref{sec:DM_definition} that this leads to results that are much more meaningful for our biological collaborators than the cP distances themselves. 
Taking into account not only the (approximative) cP distances but also the approximative correspondence maps will lead to even better distances, via the Horizontal Base Diffusion Distance. 

\begin{definition} [Horizontal Base Diffusion Distance (HBDD)]
Given a collection of surfaces $\left(S_k\right)_{k=1,\ldots,K}$, endowed with a distance between the surfaces and pairwise correspondence maps 
between the surfaces. 
Let $\left(\mathbf{W}_t\right)_{t > 0}$ be the bundle diffusion matrices constructed by means of these distances and correspondence maps, as defined in 
\eqref{Diffusion_matrix}, and let 
$\Psi_1$, $\ldots, \,\Psi_L$ be the first few dominant non-trivial eigenvectors of $\mathbf{W}_{\tau}$, for one appropriately chosen $\tau >0$. 
with corresponding eigenvalues 
$\lambda_1$, $\ldots,\, \lambda_L$. Each of these eigenvectors has an entry for each surface $k$ and for each point $x \in S_k$. For each surface $k$, we consider the matrix $G_k^{\gamma}$ of the restrictions to $S_k$ of the $L$ eigenvectors, weighted by the 
$\gamma$-th
power of their corresponding eigenvalues. 
Then the Horizontal-based diffusion distance (HBDM) 
between surfaces $k$ and $k'$ is given as a matrix-distance between $G_k^{\gamma}$ and $G_{k'}^{\gamma}$,

\begin{equation}
  \label{eq:HBDM_v1}
  \begin{aligned}
    d_{\mathrm{HBDD},\gamma}&\left(S_k,S_{k'}) \right) = \left\|G_k^{\gamma}-G_{k'}^{\gamma} \right\|.
  \end{aligned}
\end{equation}
\end{definition}

\subsection{Quantifying Errors of Shape Distances}
\label{sec:Quantifying Errors}
One of the pillar-stones of good research is the ability to measure the error of the proposed method, providing insight into the reliability and accuracy of measurements and algorithms. Error evaluation helps not only to identify weaknesses in our methods but also allows us to refine and improve them for future use. Ultimately, by rigorously evaluating errors, we can increase confidence in our findings and advance knowledge and innovation across various domains. 
Although we shall not go into detail of the error analysis for many of the results reported here, we do want to list a few of the measurement tools we use.

In the paper \cite{gao2018development} three methods are proposed to evaluate the cP registration. First, as initially proposed in \cite{boyer2011algorithms} the Mean Square Error (MSE) can be calculated as the squared Euclidean distances between the expert-determined landmarks on a surface and those propagated, via the cP correspondence map, of expert-determined landmarks on other surfaces. Larger MSEs indicate greater deviations between the computed correspondence map and the (unknown) "true"
biological correspondence map.
Second, the correspondence can be evaluated by comparing the surface characteristics, e.g. the curvature. Last  \cite{gao2018development} considered the skewness of the distribution of candidate maps, and assessed the distinctiveness of the optimal map. Although alignment errors occur with undesirable frequency (27.6\%), most of them are relatively small and represent some landmark distortion at either the anterior or posterior end of
the tooth; the most egregious errors arose when landmarks were propagated
between very dissimilar teeth. As will be shown in the next subsections, such alignment errors between pairs of surfaces with a large cP distance have very little influence on the 
next-level maps constructed by applying (horizontal) diffusion 
to the surface collection, considered as a fibre bundle.

\section{Manifold Learning}
\label{sec:Manifold Learning}
Given a data collection in which the samples are observed or assumed to vary gradually as some parameters governing the samples are varied, it is natural to view the collection as having a manifold structure. The manifold assumption provides us with the tools to study the geometry of the Riemannian manifold. Dimension reduction is a prominent approach to dealing with high-dimensional data (PCA \cite{pearson1901l}, Multidimensional Scaling \cite{cox2000multidimensional}, Linear Discriminant Analysis \cite{fisher1936use}, Locality Preserving Projections \cite{he2004locality},  Locally Linear Embedding \cite{roweis2000nonlinear}, ISOMAP \cite{tenenbaum2000global}, to mention just a few); other methods reconstruct and study the manifold geometry in high dimension (for a detailed survey about such methods please see \cite{faigenbaum2020manifold}). We focus here on Diffusion Maps \cite{coifman2005geometric}, which play a big role in our geometric analysis. In this section, we review the mathematical apparatus that is used in the next section.
Our main vehicle for studying the geometry of the base manifold is the Diffusion Maps on Fibre Bundles \cite{gao2021diffusion}, and utilize the Procrustes distance between surfaces defined in Section \ref{sec:Distance} as "input".

\subsection{Diffusion Map}
\label{sec:DM_definition}
The Diffusion Map (DM) \cite{coifman2005geometric} is a frequently used non-linear method to perform dimensional reduction. It proposes a
probabilistic interpretation for graph-Laplacian-based dimensionality reduction algorithms. Under the assumption that the discrete graph is appropriately sampled from a smooth manifold, it assigns transition probabilities from a vertex to each of its neighbors (vertices connected to it) according to the edge weights, thus defining a graph random walk, the continuous limit of which is a diffusion process over the underlying manifold. The eigenvalues and eigenvectors of the graph Laplacian,
which converge to those of the manifold Laplacian under mild assumptions
\cite{belkin2008towards, belkin2006convergence}, then reveal intrinsic information about the smooth manifold. The associated operators form a semi-group, as discussed below.


\begin{definition}[Diffusion Map]
Let $\{x_i\}$ be a data set of points in $\R^n$ sampled from a Riemannian manifold $M$ whose intrinsic dimension is low.
We define the diffusion operator on the discrete data samples as $D_t^{-1}W_t$ where
$W_{t;i,j} = \exp(|| x_{i} - x_{j} || ^2/4t) $, 
and $D_t = [D_{t;i,i}]$ is the diagonal matrix with entries defined by $D_{t;i,i} = \sum_j W_{t;i,j}$ (non-diagonal entries of $D_t$ are 0).
\end{definition}
Explicitly the diffusion map is given as
\begin{equation} H_{t;i,j} = N_{t;i} 
\exp\left({\frac{|| x_{i} - x_{j} ||}{4t}  ^2}\right)~,
\label{Diffusion_matrix_v0}
\end{equation}
where $N_{t;i}$ is a normalizing factor ensuring that 
$\sum_j \, H_{t;i,j}\,=\, 1 ~$
for all $i$. 

We omit here a detailed discussion of the theory underlying the use of the diffusion map; see e.g. \cite{coifman2005geometric, singercoifman}. When the $x_i$ correspond to a uniform sampling of points on a compact Riemannian manifold $\mathcal{M}$, the eigenvalues and eigenvectors of the matrix
$H$ approximate, in the limit as the sampling density goes to $\infty$, 
and for sufficiently small $t$, the eigenvalues and eigenvectors of the (heat) diffusion on $\mathcal{M}$ induced by the Laplace-Beltrami operator on 
$\mathcal{M}$. The argument uses that the heat diffusion kernel has the same form on all Euclidean spaces (only its normalization factor,  reflects the dimension). For $x_i$ and $x_j$ sufficiently close together on $\mathcal{M}$, the difference between the distance $\|x_i-x_j \|$ in the ambient high-dimensional Euclidean space and the distance to the origin of the pull-back to the tangent space at $x_i$ of the point $x_j$ becomes negligible, so that the kernel in (\ref{Diffusion_matrix_v0}) is appropriate in this limit. 
The main advantage of the diffusion map is that it allows to perform a local, non-linear dimensional reduction, and enables dealing with manifolds with various geometries. In practice, one carries out a spectral analysis of $H$, and one parametrizes each  $N-$dimensional data point $x_i$ 
by the $k-$tuple 
$\left( \lambda_{\ell}^{\tau}\,u_{\ell}(i) \right)_{\ell=1, \ldots,k}~$, where the $u_{\ell}$ and $\lambda_{\ell}$ are the eigenvectors and eigenvalues of $H$, $k$ is a (smallish) integer, and $\tau$ is set to a value dependent on the problem. 

\textbf{Diffusion Semi-Group Property}

\label{sec:semi-group}
Tuning the diffusion time is crucial, and is tightly connected to the question of studying local and global geometry. Picking a too large diffusion time can result in misreading local manifold geometry, since it may give too much weight to large distances, the estimates of which are expected to be much less trustworthy. On the other hand, if the diffusion time is too small, then there is no or barely any diffusion, so that the analysis fails to capture the geometric properties of the underlying manifold. In order to choose the optimal diffusion time, the paper \cite{shan2022diffusion, shan2019probabilistic} proposed to use the extent to which the candidate-diffusion operator exhibits Semi-group behavior in order to calibrate the $t$. Let us first define the Semi-group Property as

\begin{definition}[Semi-group Property]
\label{semi-group property_def}
    A family of operators $\{T^t\}$ for $0\leq t< \infty$, in which $T^0$ is the identity operation,  is said to have semi-group property if $T^{t_1+t_2} = T^{t_1}T^{t_2}$.
\end{definition}

The family of heat diffusion operators on a manifold forms a strongly continuous semigroup. The discrete operators $H_t$ can be reasonably expected to be good discretizations of a true
manifold diffusion operator only if they approximately obey the Semi-Group Property. This motivates the following definition:

\begin{definition}  [Semi-group Error]
  The ($t$-dependent) Semi-group Error of the operators 
  $\left(H_t\right)_{t>0}$ is defined as
  \[
  SGE(t) = ||(H_{t}^2 - H_{2t}||
  \]

\end{definition}

If $\left(H_t\right)_{t>0}$ were a true semi-group, then $SGE(t)$ would be zero for all $t$. 
As shown in \cite{shan2022diffusion, shan2019probabilistic},
in practice 
$SGE(t)$ typically decreases as $t$ increases (starting from $t$ close to zero), then reaches a minimum, after which it starts increasing again. Picking the diffusion time $t$ around the value for which this minimum is reached is a heuristic way of "tuning" the diffusion time that produces good results in applications, as illustrated in \cite{shan2022diffusion, shan2019probabilistic}.  In what follows, we shall assume that
the value $t$ picked for a diffusion operator $H_t$ 
has been "tuned" in this way.

\subsection{Diffusion Distance}
\label{sec:DM_definition}
Although the operator $H_t$ is not symmetric (because of its non-symmetric normalization), one can determine its
eigenvectors and eigenvalues by simply considering the 
symmetric operator $D_t^{1/2} H_t D_t{-1/2}$; this symmetric operator has an orthonormal basis of eigenvectors, applying $D_t{-1/2}$ to each of these eigenvectors turns it into an eigenvector of $H_t$, with the same eigenvalue. Typically only the first non-trivial dominant eigenvectors of $H_t$ are 
"trustworthy": the distance matrix estimate from which we start is typically only an approximate estimate, and their
inaccuracy affects the higher-order eigenvectors and eigenvalues more than the first few. 

Note that we use only the non-trivial eigenvectors: because 
of the normalization of $H_t$, its very first eigenvector is the
constant-1 vector, with eigenvalue 1; it is non-informative
and typically not considered. Given those eigenvectors and eigenvalues, the diffusion distance is later computed as indicated in \eref{eq:DD}. The non-linear embedding identifies a new domain that represents meaningful geometric descriptions of data sets, by generating efficient representations of complex geometric structures.

\begin{figure}[h]
  \centering
  \includegraphics[width=.35\textwidth]{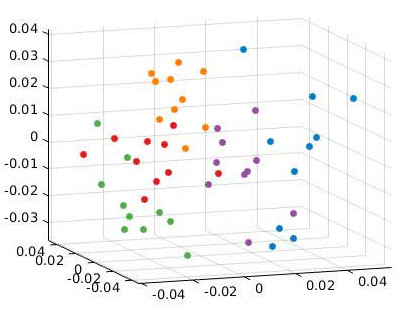}
  \includegraphics[width=.4\textwidth]{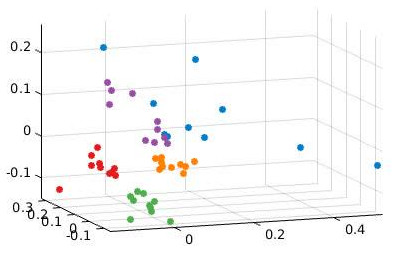}
 
  \caption{Comparing the cP distances (left) and DD distances (right)
for a dataset of second 
molars from 50 species belonging to 5 genera of monkeys; in each case a 3D MDS representation of the 
distances are shown, with one colored dot for each molar, where
species belonging to the same genus are given the same color.}
  \label{DD_vs_cP}
\end{figure}

This can be seen in Figure \ref{DD_vs_cP} which illustrates the power of this diffusion approach. Its left panel shows a multi-dimensional 
scaling (MDS) representation in 3D of the cP distances for a data set of 50 second
molars of different monkey species, belonging to 5 genera (10 species per genus). The right panel of the figure uses those
{\em same} cP distances as input, but analyzes them via the manifold diffusion processing described above (and in greater detail in \cite{gao2015hypoelliptic}); the figure shows the MDS representation in 3D of the corresponding Diffusion Distances computed from this input. It is clear that this extra step in the analysis provides a much better separation of the genera. 

\vspace{-0.0in}

\subsection{Manifolds in which each Sample has Structure: Fiber Bundles}
\label{sec:manifold learning_Fiber}
The previous subsection illustrated that diffusion distances provide a powerful tool to analyze (dis)similarities of teeth and bones in a biological morphology framework. In using the
cP distances as input for a manifold diffusion analysis, we 
took into account only part of the information produced by
the continuous Procrustes approach: as indicated by the formula \eqref{eq:cont_P_dist}, the cP distance between surfaces $S_1$ and $S_2$ is the minimal value obtained in an optimization procedure that also
determines a corresponding minimizing correspondence map 
between $S_1$ and $S_2$. In the Horizontal Base Diffusion Distance (HBDD) construction, we use not only the cP distances
but also the correspondence amps as input; we do this by considering the unknown manifold describing (or
parametrizing, or underlying) the collection
of surfaces as a {\em base manifold}, and each of the surfaces themselves as the {\em fibers} in a {\em fiber bundle} framework. In the next few paragraphs, we briefly review this
framework, in which we can naturally define a diffusion analysis that lends itself to our purposes.

\begin{definition}[Fibre Bundle] Assume that 
$M$ and $F$ are orientable Riemannian manifolds. A \emph{Fibre bundle} with \emph{base space} $M$ and \emph{fibre space} $F$ is a quadruple ($E, M, F, \pi$) where $E$ (called the \emph{total space}) is a differentiable Riemannian manifold and $\pi : E \rightarrow M$ is a surjective smooth map from $E$ to $M$ that satisfies the local triviality condition, meaning that for every $p\in E, \exists$ an open neighborhood $U\subset M$ of $\pi(x)$ and a diffeomorphism \[\phi: \pi^{-1}(U) \rightarrow U\times F\]
such that $\pi|_U = \pi_1 \circ \phi$, where $\pi|_U$ is the restriction of $\pi$ to $U$ and $\pi_1$ (from $U\times F$ to $U$) is the projection to the first components (i.e. it maps $(u,f)$ to $u$). 

\end{definition}

For any point $p\in M$, $\pi^{-1}(\{p\})$ is diffeomorphic to $F$. We then denote $F_p =\pi^{-1}(\{p\})$ as the \emph{fibre} over $p$. The local trivialization allows a way to describe a specific fibre in a space with a simpler topology. 

In our case, we also assume that the bundle comes equipped with correspondence maps between fibres, $P_{yx}: F_x\rightarrow F_y$. This function maps  points on $F_x$ to points on $F_y$ in a smooth, bijective way. It associates the two fibres and identifies their points to each other. (In \cite{shan2019probabilistic} these are interpreted as \emph{parallel transport} in a probabilistic setting.) 

\begin{figure}[h]

 \centering
	\captionsetup[subfloat]{farskip=0pt,captionskip=0pt, aboveskip=0pt}
	\subfloat[][]{ \includegraphics[width=0.4\textwidth]{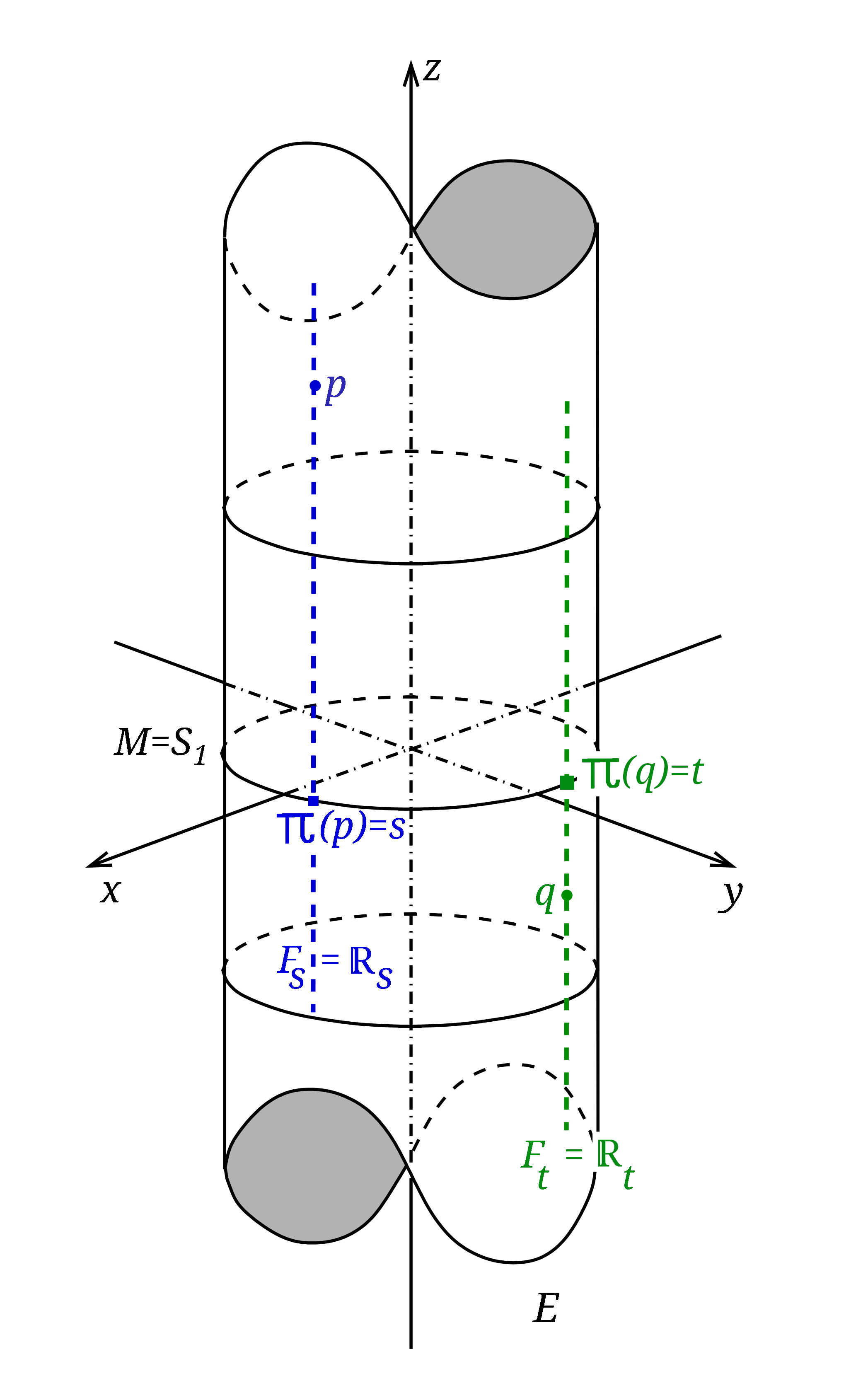} } 
 \hspace{2em}
	\subfloat[][]{ \includegraphics[width=0.45\textwidth]{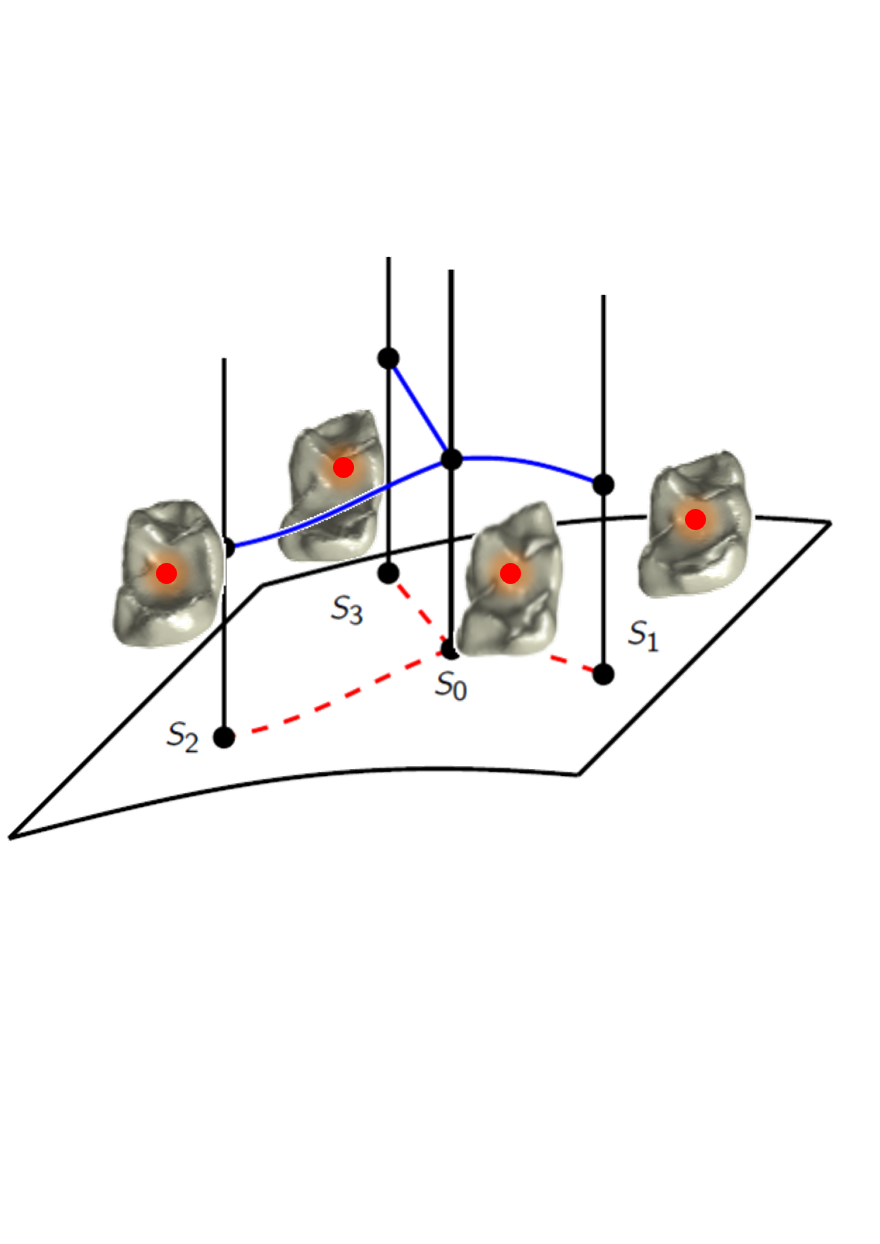} } 
 
	\caption{Left: Illustration of the decomposition of the cylinder manifold into the base manifold $S_1$, and the vertical fibers, where the total space can be represented as a Cartesian product of the base manifold, and the fiber; Right: Fiber bundle of teeth. Each anatomical surface corresponds to a point on the base manifold. Each surface constitutes a fibre.}.
	\label{fig:Fibres}
\end{figure}

An illustration of these concepts appears in the toy example illustrated in the left panel of Figure \ref{fig:Fibres}. Here  the two-dimensional manifold $E=\{(x,y,z) \in \R^3; \,x^2+y^2=1\}$ is a cylinder. Then one can view the circle $S_1$ 
(the cross-section of the cylinder) can be viewed as the base manifold, and $\R$ is the fibre. The coordinates $x,y,z$ of $p\in E$ can be viewed as a combination of $s:=(x,y) \in S_1$, $z \in \R$. Now, the projection map $\pi: E \to S_1$ is defined as 
$\pi (x,y,z) = (x,y) = s$ and its inverse
$\pi(x)^{-1} \cong\R_s$ i.e. the entire line (the fibre over $s\in S_1$). (Note: in this case, $E$ can be identified simply with the product space of $M$ and $F$; in cases like these, the bundle is said to be {\em trivial}. For general fibre bundles this need not be the case -- for instance, one can add a slight twist to this toy construction and obtain a Moebius-like total space $E$ that is no longer diffeomorphic to 
$M \times F$.)

In our setting, the total space $E$ consists of the collection of all the points in the union of a collection of (tooth) surfaces that we assume can be parametrized smoothly by a base manifold $M$ that we want to study. As long as we look at species that are closely related to each other, we believe that we can bring 
them back, following a reasonably short path in a phylogenetic  tree, to a common ancestor with reasonably similar teeth. In that case, the mapping from one surface $S_i$ to another $S_j$ is the same as the mapping from $S_i$ to an ancestral surface 
$S_A$, and back from $S_A$ to $S_j$. This would then mean that the mappings are all consistent, over the whole collection -- that mapping from $S_i$ to $S_j$ and then from $S_j$ to $S_k$ produces the same result as directly from $S_i$ to $S_k$. (When we will start looking at surfaces for more distantly related species, this will no longer be the case, and we expect our fibre bundles to become non-trivial.) This is illustrated in the right panel of Figure \ref{fig:Fibres}. 

\textbf{Horizontal random walks and diffusion processes on fibre bundles}

\label{HDM_definition}
We shall use the fibre bundle model to define diffusion operators on the total space $E$, and hope to learn from its spectral analysis more than from simple diffusion on the 
base manifold $M$, using the cP distances as input. 

To this end, we have to define distances on $E$, or,
equivalently, define a diffusion kernel on $E \times E$. 
In practice, we have on each of the surfaces $S_k$, with
$k$ ranging from 1 to $K$, a collection of $N_k$ 
sample points 
$\left(x_{i,k}\right)_{i=1,\ldots,N_k}$, where the 
$N_k$ for different $k$ could be (but need not be) equal to each other; for simplicity, we shall assume they are equal in this exposition. 
We also have (from the computational optimization
inspired by the theoretical definition of cP distance), for
each pair of surfaces $k,\,\ell$,
correspondence maps $P_{k,\ell}$ that map the collection
$\left(x_{i,k}\right)_{i=1,\ldots,N}$ bijectively 
to $\left(x_{j,\ell}\right)_{j=1,\ldots,N}$ (we are assuming
$N_k=N_{\ell}=N $ here); we shall consider these $P_{k,\ell}$ to be (noisy) approximations for the (unknown) correspondence maps between the continuous surfaces. "Random Walking" on the total surface will consist of walking from surface to surface and, on
each surface, from point to point. We don't have an a priori way
to "synchronize" the walking in the fibre direction with that along the base direction, so we shall use two different time constants, and explore later what it means to vary one vs. the other. Taking all this into account, we define the (discrete and approximate) diffusion kernel on $E$ as

\begin{equation} \mathbf{W}_{\tau_1,\tau_2;j,\ell,i,k} = N_{j,\ell} \exp\left(
-\frac{|| P_{k\ell}(x_{\ell, j}) - x_{k, i} ||^2 + 
|| x_{\ell,j} - P_{\ell k} (x_{k,i}) ||^2}{4\tau_1} \right) 
\exp\left(-\frac{D_{k,\ell}^2}{2\tau_2}\right) 
\label{Diffusion_matrix}
\end{equation}

where $D_{k,\ell}$ denotes the cP distance between $S_\ell$ and $S_k$ and $P_{k\ell}$ is the map from 
$\left(x_{j',\ell}\right)_{j'=1,\ldots,N}$ 
to
$\left(x_{i',k}\right)_{i'=1,\ldots,N}$ described above; the normalization 
$N_{j,\ell}$ ensures that the sum over all $i$ and $k$ of the 
$\mathbf{W}_{\tau_1,\tau_2;j,\ell,i,k}$ equals 1.
The diffusion parameters $\tau_1$ and $\tau_2$ control the amount of information that is taken into account from the base manifold versus from the fibers. 

Thus the HDM matrix is a block matrix in which the $(\ell, k)$-th block captures the comparison of surface $\ell$ and surface $k$ while relying on the correspondence map between them.
The $(j,i)$ entry in the $(\ell, k)$-th block of the HDM matrix stands for the edge weight between $x_{\ell,j}$, and $x_{k, i}$.
Choosing $\tau_1$ and $\tau_2$ can be tricky; increasing one while keeping the other fixed corresponds to diffusing more rapidly in one of the two directions (base manifold or fibre) 
than the other. We'll assume that both values are chosen in a regime where the semigroup error is small (see above); within that region, some latitude is possible. Depending on what we want to explore, we may
adjust the ratio $\tau_1/\tau_2$; we shall come back to this below.

\textbf{Spectral Distances and Embeddings}

The eigenvectors $\Psi_r$ of $\mathbf{W}$ have 
$N_1+N_2 + \ldots+ N_K$ entries 
(or $KN$ if all the $N_k$ are 
equal to $N$). 
As expected from a diffusion kernel, due to the normalization of the rows of $H$, the largest eigenvalue is a trivial eigenvector, with all entries equal to $1$, and eigenvalue $1$. 
We shall systematically disregard this (uninformative) eigenvector $\Psi_0$, and consider only $\Psi_\ell$ with $\ell \geq 1$. We can consider the eigenvectors as concatenations of $K$
segments of lengths $N_k$, each indexed according to surface 
$S_k$, with $k$ ranging from 1 to $K$; let's denote these
segments by $\psi_{\ell,k}$.(Equivalently, for $m=1, \ldots, N_k$,
$\left(\psi_{\ell,k}\right)_m = \left(\Psi_\ell\right)_{N_1+\ldots+
N_{k-1}+m}$.) Because the data from which the spectral analysis
is constructed were noisy (the points picked on the different $S_k$ are not in strict correspondence to each other, since they were picked randomly on each surface, independently for each surface, leading to errors not only in the correspondence maps but also in the approximate cP distances between the surfaces), the eigenvectors are noisy as well; this expresses itself all the more as the corresponding eigenvalue decreases. In practice, we therefore use only the $L$ dominant eigenvectors and eigenvalues, i.e. we assume the ordering is such that the $\lambda_\ell$ decrease as $\ell$ increases, and we use only the  $\Psi_\ell$ and $\lambda_\ell$, for $\ell \leq L$; the
value we pick for $L$ can vary depending on the application. Below we shall see two applications, one with $L=3$ and another with $L=100$; in all cases $L$
is significantly smaller than any of the $N_k$. 

One way to "distill" the geometric content from the $\Psi_\ell$,
$\ell=1,\ldots L$, as it pertains to each of the surfaces $S_k$,
is to consider, 
for each $k$, the $R\times R$ Gramm matrix $G_k$ defined by 
\[
\left(G_k\right)_{\ell,\ell'} = \langle \psi_{\ell,k} , \psi_{\ell',k} \rangle = \sum_{m=1}^{N_k}\, 
\left(\Psi_\ell\right)_{N_1+\ldots+
N_{k-1}+m} \overline{\left(\Psi_{\ell'}\right)_{N_1+\ldots+
N_{k-1}+m}}~~;
\]
we also consider the weighted Gramm matrices $G_k^t$ defined by
\[
\left(G_k^{\gamma}\right)_{\ell,\ell'}= \lambda_{\ell}^{\gamma/2} \lambda_{\ell'}^{\gamma/2}\left(G_k\right)_{\ell,\ell'}~~.
\]

In \cite{gao2021diffusion}, the  \emph{HBDD distance}
$d_{\mathrm{HBDD},\gamma} \left( S_k , S_{k'} \right)$ 
is then defined as a matrix distance between 
the two $R \times R$
matrices $G_k^{\gamma}$ and $G_{k'}^{\gamma}$. (Several possible candidates can be considered for this matrix distance; if one uses the Hilbert-Schmidt norm, then the HBDD distance can be easily reformulated in terms of the original  $\psi_{\ell,k}$ vectors, which is then closer to the embeddings considered in \cite{coifman2006diffusion}, extended to this fiber bundle setting.) 

\begin{figure}[h]
  \centering
  \includegraphics[width=.33\textwidth]{images/cPDMDS.jpg}
  \includegraphics[width=.35\textwidth]{images/DMMDS.jpg}
  \includegraphics[width=.31\textwidth]{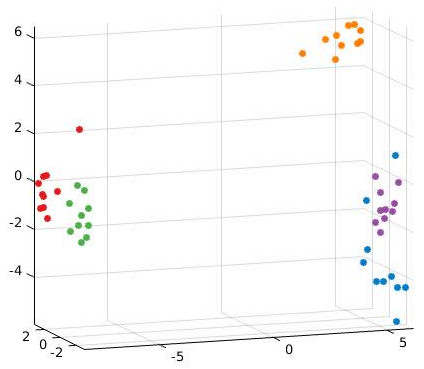}
  \caption{Geometric Similarity of Anatomical Surfaces using different distance calculations.   Left: Embeddings into $\mathbb{R}^3$ using \emph{Multi-dimensional Scaling} of CPD  \eqref{eq:cont_P_dist}; Middle: DM calculated after the diffusion map embedding \eqref{Diffusion_matrix_v0}; Right: and HBDD.}.
  \label{fig:MDS_comparison}
\end{figure}

To illustrate the improvement of the HBDD distance over the Diffusion Distance, we revisit the comparison provided in Figure \ref{DD_vs_cP}, adding now a third panel on the right in which we show the MDS embedding of the HBDD distances for the same data set. It is apparent that the surfaces (2nd molars of monkeys in this case) for the five different genuses are clustered much better with HBDD than with the Diffusion distance. In addition, the HBDD reflects the dietary groups within the data set: the folivores Alouatta (red) and Brachyteles (green) are adjacent to each other in the rightmost panel of Figure~\ref{fig:MDS_comparison}, and so are the frugivores Ateles (blue) and Callicebus (purple); the insectivore Saimiri (orange) is far from the other herbivorous groups.

\subsection{Improved Spatial Correspondence}
\begin{figure}[h]
	\centering
	\label{fig:c}\includegraphics[width=1\textwidth,keepaspectratio]{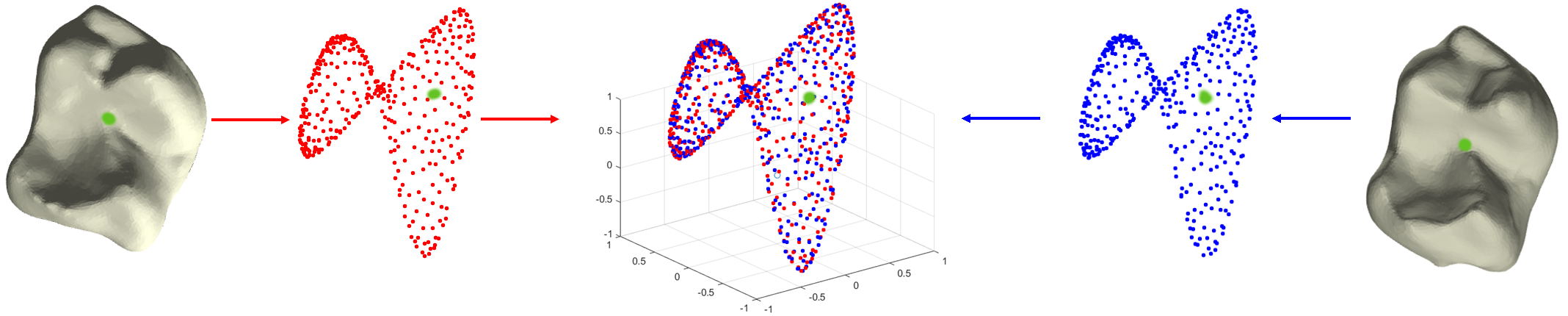} 
	
	\caption{Fiber bundle diffusion map calculated based on a collection of teeth. Its eigenvectors are later used to embed the surface coordinates into 3-dimensional space (in the middle). The different colors (red, and blue) represent points from the two different teeth projected to the diffusion space. As can be seen, the projected surfaces closely recover the fiber manifold. }
	\label{fig:teeth_example_pringle}
\end{figure}

 The spectral analysis of the Bundle Diffusion provides a different parametrization for the surfaces: for each $k=1,\ldots,K$, we associate to each point 
 $p_m \in S_k$, with $m=1,\ldots,N_k$ the $L$-tuple 
 $\left(\psi_{\ell,k}\right)_m = \left(\Psi_\ell\right)_{N_1+\ldots+
N_{k-1}+m}$, with $\ell$ ranging from 1 to some predetermined $L$. This defines a mapping from the entire collection of surfaces into $\mathbb{R}^L$; we call this the Horizontal Diffusion Map (HDM). See Figure \ref{fig:teeth_example_pringle} for an illustration of embedding two surfaces (outmost left and right) into the diffusion space to recover the fiber manifold (in the middle). The main advantage of HDM embedding is that not only does it uncover a new domain that captures meaningful geometric characteristics of datasets, but also preserves the intrinsic geometry of each surface $S_k$.

In this section, we review a few different applications of the HDM,  related to improving spatial correspondences between the $S_k$.

\textbf{Slice and Dice: Segmenting Surfaces to Match Evolutionary Equivalents}
\label{sec:Segmenting Surfaces}

\begin{figure}[h]
  \centering
  \includegraphics[width=.9\textwidth]{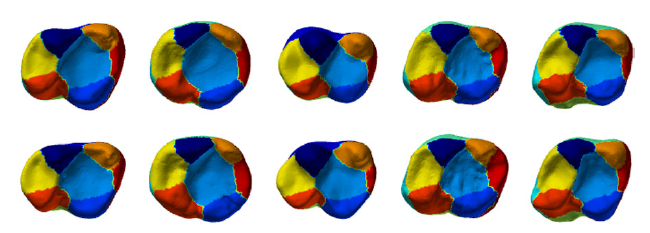}
  \caption{Globally consistent segmentation by spectral clustering for HDM. From the left to the right row: Alouatta, Ateles, Brachyteles, Callicebus, Saimiri, image source \cite{gao2021diffusion}. In each column are different teeth examples from the same species.}
  \label{fig:hdm_spectral_clustering}
\end{figure}

Once we have embedded a collection of surfaces into $\mathbb{R}^L$, we can use standard clustering algorithms to determine the best way to separate the collection into some predetermined number $M$ of clusters.  
Figure~\ref{fig:hdm_spectral_clustering} shows the result
on 10 different primate teeth (part of a larger collection of
32 that were analyzed with Bundle Diffusion); in this case 
$L$ was picked to be 100, and $M=12$. To visualize the clustering done in HDM space in terms of the surfaces 
themselves, one of 12 colors was assigned to each of the 12 clusters, and points on the surfaces in the collection were given the color of the cluster to which their HDM embedding 
belonged. 

\vspace{-0.1in}

Although there is no a priori constraint imposing it, each cluster color is represented on each of the $S_k$, and the zones that they delimit on the different teeth appear to be in
biological homologous correspondence. Note that the same-color zones on different teeth do not have the same area -- even though the analysis started from correspondence maps 
(the cP maps) that were (approximately) area-preserving. 
Our biological colleagues are very excited by this potential of the HDM method to segment the original meshes in a consistent way across the entire shape space, similar
to what is done in classical anatomical studies, where features are commonly not addressed as points (or landmarks), but rather as sub-regions of the anatomy, and they are starting to use this as a tool in their own statistical analyses. 

\textbf{Improved Spatial Registration} 

A different but likewise very important application uses
a very low value of $L$, more precisely $L=3$. Upon 
examination of the resulting embedding into $\mathbb{R}^3$, 
one finds that all the $S_k$ are mapped to essentially the
same 2-dimensional surface
(up to noise effects), as illustrated in Figure \ref{fig:teeth_example_pringle}. We call this surface the \textit{template surface}. 

Of course, the images on the template surface of the randomly picked points $p_{k,m}$ on each $S_k$ are typically all
different. By using the template surface, it is now possible, however, to compute for each point $p_{k,m}$ on  $S_k$ a corresponding interpolation of 
the points $p_{k',m'}$ on another $S_k'$ that is much better
correspondence with $p_{k,m}$ than any of the individual $p_{k',m'}$. To do this, one considers the image $P_{k,m}$
on the template
space, under the embedding, of $p_{k,m}$, and determines weights for the nearby images $P_{k',m'}$ that would make 
their weighted average coincide with $P_{k,m}$. Taking the same
weighted average of the $p_{k',m'}$ on the original $S_{k'}$ 
defines then a point on $S_{k'}$ that is a better candidate for
the point on $S_{k'}$ "truly" corresponding to $p_{k,m}$. 
One can then use these improved correspondence maps to repeat
the entire procedure of constructing the Bundle Diffusion,
identifying its Eigenvectors, determining a template surface, and going through another step of improving the correspondence maps.
Numerical experiments on toy examples suggest that it may be
useful to repeat this iterative procedure a few times; we
don't know whether or under what conditions it would converge.

\textbf{Collection-Based Gaussian Process Landmarking}

In subsection \ref{sub_sec_regi_align} we described Gaussian-Process Landmarking, a technique that 
identifies "geometrical landmarks" on individual surfaces by selecting successively those points that contribute the most new information about the surface, given the points chosen earlier. 
Once the surfaces $S_k$ in a collection have been jointly registered
via the HDM (as described above), one can use a similar approach for the whole collection. In this case one defines a Gaussian kernel on the template surface using the covariance between points over {\em all} the surfaces $S_k$ simultaneously, and then selects
successive landmark points on the template surface, each time picking a new landmark as the point with the largest variance conditioned on all the previously selected
landmarks. This automatically defines corresponding landmarks
on all the $S_k$ via the HDM map. This is described in more detail
in \cite{Gaussian_Process_Joint_Landmarking_of_Surfaces}.

\section{Other differential geometry tools}
\label{sec:ariaDNE}

The preceding sections have concentrated on tools from 
differential geometry that we brought to bear on the analysis
{\em collections} of surfaces. However, we have also developed
tools to reliably compute geometric properties for individual
surfaces; in fact, the very first collaboration of our team combining applied mathematicians and biological morphologists \cite{bunn2011comparing,winchester2014dental,
winchester2016morphotester, pampush2016introducing} concerned the use of 
effective discrete algorithms to compute the Dirichlet Normal Energy (DNE) of a surface and its use on molar shapes
to distinguish dietary
preferences of different mammal species. 
The DNE is a measure of curvature, computed as 
the integral, over the 
surface (normalized to have area 1), of the change in direction 
of the normal vector. Compared with some other methods used in 
morphological studies, the DNE had the 
advantage of being landmark-free and independent of the surface’s initial position, orientation,
and scale, resulting in a more robust metric, less sensitive to noise than some other tools. 
Nevertheless, our first implementation for the calculation of the DNE was strongly tied to the mesh representing the surface; 
once it was applied to a wider range of mesh preparation protocols, it turned out to be overly sensitive to variation
in some of the mesh preparation
parameters, especially for the boundary
triangles. 

In \cite{shan2019ariadne} a more robust implementation 
of the same mathematical concept was introduced, called {\em ariaDNE} (for {\em {\bf a} {\bf r}obustly {\bf i}mplemented} {\em {\bf a}lgorithm} for {\em {\bf D}irichlet
{\bf N}ormal {\bf E}nergy}). 

To compute the DNE numerically, ariaDNE uses a local PCA to determine the tangent plane as well as the curvature (both
approximately) --  at sample points $P$ where the surface 
does not have 
overly sharp ridges or cusps, the first two principal components of the PCA of the sample points in a neighborhood of $P$ are tangent to the surface, the third principal component approximates the normal direction, and the smallest principal value corresponds to the curvature. Near a ridge or
cusp, it may be necessary to adjust the procedure, since the
normal direction may then no longer be the one with the smallest PCA
eigenvalue; in this case, it is however easy to determine
which of the three eigenvectors is normal to the surface
at $P$, and select the corresponding eigenvalue. 
For numerical stability, ariaDNE uses weighted PCA (with
a Gaussian kernel centered on $P$) rather then local
PCA. 

\begin{figure}[h]
  \centering
  \includegraphics[width=.8\textwidth]{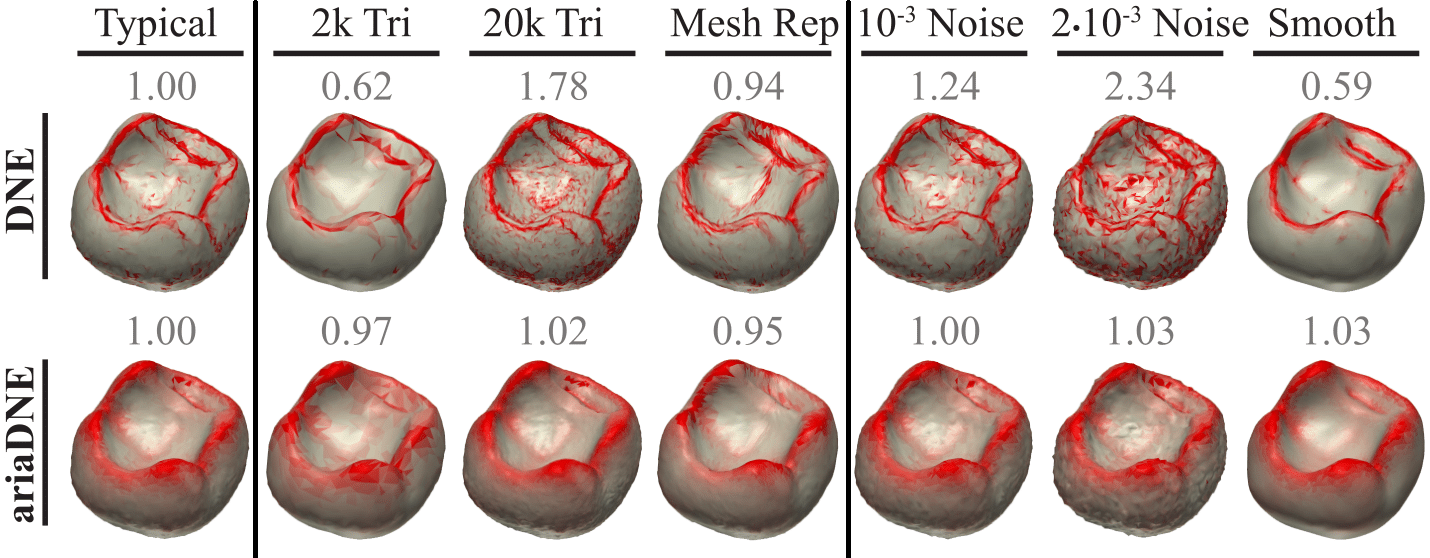}
  \caption{Comparing effects of triangle count, mesh representation, noise, and
smoothing on  DNE (top) and ariaDNE (bottom). The distribution of curvature
as measured by each method overlaid in shades of red on a grey 3D rendering
of the surface is presented \cite{shan2019ariadne}. }
  \label{fig:form_function}
\end{figure}

Figure \ref{fig:form_function} shows DNE and ariaDNE values on an example tooth and various meshes representing the same
surface and different amounts of noise: (from left to right) 2k triangles, 20k triangles, a different mesh representation,
0.001 noise magnitude, 0.002 noise, and smoothing (commonly used to eliminate noise produced during mesh
preparation). The red surface shading indicates
the value of curvature as measured by each method.  As can be seen, by the max curvatures value heading the surfaces, the ariaDNE is remarkably stable to noise and different mesh configurations.

\vspace{-0.1in}

\section{Summary and Future Research Directions}

This paper summarizes methods that were developed, based on differential geometry and machine learning, to analyze collections of anatomical surfaces within the framework of evolutionary anthropology. 
The surfaces, representing teeth or bones, are scanned and processed as triangular meshes; they pose challenges, due to nuanced anatomical variations, that are quite different from the classical problems in Computer Graphics and CAGD. To address these
challenges, tools are required that differ significantly from classical surface analysis methods. Through a comprehensive survey spanning more than 15 years of research, the paper covers methods for surface comparison, registration, alignment, and distance calculation. It also addresses broader topics such as landmark selection, shape segmentation, and classification analysis, and showcases the advantages of manifold analysis in addressing evolutionary biology inquiries.
While focusing primarily on disc-type surfaces resembling teeth, the methods discussed are adaptable (and in several instances, have been adapted) to sphere-type surfaces like bones, suggesting broader applicability and potential avenues for future research. We are currently working also on combining some of the approaches outlined above, e.g. in joint Gaussian Landmarking for collection
of surfaces and iterative increasingly accurate registration (both mentioned above), but also incorporating not only
spatial coordinates but also other features (such as the scalar DNE field) in computing generalized continuous Procrustes distances
to be used in subsequent diffusion-based learning. We also hope to
learn how to coarse-grain the analysis so as to extend our
analysis to species that are more distantly related. Finally, we are convinced our bundle-diffusion approach is applicable to many other manifold learning settings beyond anatomical surfaces, in which the data "points" have internal structure and the (dis)similarity between them relies on comparing these internal structures.

The tools outlined in this paper are summarized in the following table (Table \ref{tab:code}), providing references to both the original paper and a link to the corresponding code repository.

\begin{table}[H]
\caption{Studying Morphological Variation: Available Code}
\label{tab:code}   
\begin{tabular}{p{1.5cm}|p{4.8cm}|p{2cm}|p{6.5cm}}
\hline\hline\noalign{\smallskip}
Method & Short description & Read more & Link to the code  \\
\noalign{\smallskip}\svhline\noalign{\smallskip}

cPD & Find the best correspondence of a pair of surfaces via optimization process &  Sec. \ref{sec:Distance},  \cite{al2013continuous} & \url{https://github.com/trgao10/cPdist}\\
\hline
 
Auto3dgm & Automated 3D Geometric Morphometric Aligment & Sec. \ref{sec:reg_landmark},  \cite{boyer2015new}&   \url{https://github.com/ToothAndClaw/Auto3dgm\_Python} \\
\hline
Gaussian Process Landmarking & Landmarking via Gaussian Kriging method, and later registration via landmark matching & Sec. \ref{sub_sec_regi_align},  \cite{gao2019gaussian} & \url{https://github.com/shaharkov/GPLmkBDMatch}\\
\hline
HDM  & Horizontal random walk on the fibre bundle & Sec. \ref{HDM_definition},  \cite{gao2021diffusion} & \url{https://github.com/trgao10/HDM}\\
\hline
ariaDNE & Robust curvature estimation  &  Sec. 
\ref{sec:ariaDNE},  \cite{shan2019ariadne} & \url{https://sshanshans.github.io/articles/ariadne.html}\\
\hline
Hecate & Consistent segmentation of biological surface regions by spectral clustering from a hypoelliptic diffusion map, and Continues Procrustes distance, and HDM  & Sec. \ref{sec:Segmenting Surfaces},  \cite{fulwood2021reconstructing} & \url{https://github.com/JulieWinchester/hecate}\\
\hline
Eyes on the Prize (EOP) Registration - SAMS & Matlab package to automatically process, register, and analyze collections of anatomical surfaces at the resolution higher   & Sec. \ref{sub_sec_regi_align},  \cite{ravier2018eyes, ravier2018algorithms} & \url{https://github.com/RRavier/SAMS}\\
\hline
SlicerMorph & Platform to retrieve,
visualize, and analyze 3D morphology
 &  \cite{rolfe2021slicermorph}   & 
  \makecell[l]{\url{https://slicermorph.github.io/\#six} \\ \url{https://toothandclaw.github.io/installations/}} \\

\noalign{\smallskip}\hline\hline\noalign{\smallskip}
\end{tabular}

\end{table}
While most of the 3D surfaces in this paper appear in the MorphoSource \cite{morphosource}, there are other datasets that we thought are worth mentioning.

\begin{table}[H]
\caption{Studying Morphological Variation: 3D specimen data}
\label{tab:data}       
%
%

\begin{tabular}{p{2.3cm}|p{4cm}|p{0.8cm}|p{7.3cm}}
\hline\noalign{\smallskip}
Data Set name & Short description & Ref. & Link \\
\noalign{\smallskip}\svhline\noalign{\smallskip}

MorphoSource & 3D meshes, and CT-images of biological skeletal material& \cite{boyer2016morphosource} & 
\makecell[l]{ \url{https://www.morphosource.org/users/sign\_in?locale=en}\\
visualization: \\
 \url{https://gaotingran.com/research/auto3dgm.html}}

\\ \hline
Second mandibular surfaces &  116 Second mandibular molars of prosimian primates and non-primate, as well as observer-determined landmark & \cite{boyer2011algorithms} &  \url{https://www.wisdom.weizmann.ac.il/~ylipman/CPsurfcomp/}\\

\noalign{\smallskip}\hline\noalign{\smallskip}
\end{tabular}

\end{table}


\section{Acknowledgements}
The authors gratefully acknowledge support from the Simons Foundation
Math + X grant 400837. Shira is grateful to the Eric and Wendy Schmidt Fund for Strategic Innovation, the
Zuckerman-CHE STEM Program.
\input{references}

\end{document}

%% file: references.tex
%
%
%

%% file: main.bbl
\begin{thebibliography}{99.}%
\bibitem{yapuncich2019digital}Yapuncich, G., Kemp, A., Griffith, D., Gladman, J., Ehmke, E. \& Boyer, D. A digital collection of rare and endangered lemurs and other primates from the Duke Lemur Center. {\em PloS One}. \textbf{14}, e0219411 (2019)



\bibitem{boyer2008relief} Boyer, D.: Relief index of second mandibular molars is a correlate of diet among prosimian primates and other euarchontan mammals. {\em Journal Of Human Evolution}. \textbf{55}, 1118-1137 (2008)



\bibitem{simons1961dentition}Simons, E. The dentition of Ourayia: its bearing on relationships of omomyid prosimians.  (1961)
\bibitem{rose1984gradual}Rose, K. \& Bown, T. Gradual phyletic evolution at the generic level in early Eocene omomyid primates. {\em Nature}. \textbf{309}, 250-252 (1984)
\bibitem{rasmussen1995dentition}Rasmussen, D., Shekelle, M., Walsh, S. \& Riney, B. The dentition of Dyseolemur, and comments on the use of the anterior teeth in primate systematics. {\em Journal Of Human Evolution}. \textbf{29}, 301-320 (1995)

\bibitem{morphosource} Morphosource, \url{https://www.morphosource.org/users/sign\_in?locale=en}, March 27, 2024.

\bibitem{GPLmkBDMatch} GPLmkBDMatch, \url{https://github.com/shaharkov/GPLmkBDMatch}, March 27, 2024.



\bibitem{boyer2016morphosource}Boyer, D., Gunnell, G., Kaufman, S. \& McGeary, T. Morphosource: archiving and sharing 3-D digital specimen data. {\em The Paleontological Society Papers}. \textbf{22} pp. 157-181 (2016)


\bibitem{Avizo} Avizo, \url{ www.thermofisher.com}, March 27, 2024.

\bibitem{dragonfly} dragonfly, \url{https://www.theobjects.com/dragonfly}, March 27, 2024.

\bibitem{nagai2015tomographic} Nagai, Y., Ohtake, Y. \& Suzuki, H. Tomographic surface reconstruction from point cloud. {Computers and Graphics}. \textbf{46} pp. 55-63 (2015)

\bibitem{boyer2015new}Boyer, D., Puente, J., Gladman, J., Glynn, C., Mukherjee, S., Yapuncich, G. \& Daubechies, I. A new fully automated approach for aligning and comparing shapes. {\em The Anatomical Record}. \textbf{298}, 249-276 (2015)

\bibitem{lipman2011conformal}Lipman, Y. \& Daubechies, I. Conformal Wasserstein distances: Comparing surfaces in polynomial time. {\em Advances In Mathematics}. \textbf{227}, 1047-1077 (2011)

\bibitem{boyer2011algorithms}Boyer, D., Lipman, Y., St. Clair, E., Puente, J., Patel, B., Funkhouser, T., Jernvall, J. \& Daubechies, I. Algorithms to automatically quantify the geometric similarity of anatomical surfaces. {\em Proceedings Of The National Academy Of Sciences}. \textbf{108}, 18221-18226 (2011)


\bibitem{al2013continuous}Al-Aifari, R., Daubechies, I. \& Lipman, Y. Continuous Procrustes distance between two surfaces. {\em Communications On Pure And Applied Mathematics}. \textbf{66}, 934-964 (2013)

\bibitem{Auto3dgm} Auto3dgm, \url{https://github.com/ToothAndClaw/Auto3dgm_Python}, March 27, 2024.

\bibitem{gao2019gaussian}Gao, T., Kovalsky, S., Boyer, D. \& Daubechies, I. Gaussian process landmarking for three-dimensional geometric morphometrics. {\em SIAM Journal On Mathematics Of Data Science}. \textbf{1}, 237-267 (2019)

\bibitem{gao2019gaussian_2}Gao, T., Kovalsky, S. \& Daubechies, I. Gaussian process landmarking on manifolds. {\em SIAM Journal On Mathematics Of Data Science}. \textbf{1}, 208-236 (2019)

\bibitem{turner2014persistent}Turner, K., Mukherjee, S. \& Boyer, D. Persistent homology transform for modeling shapes and surfaces. {\em Information And Inference: A Journal Of The IMA}. \textbf{3}, 310-344 (2014)


\bibitem{moenning2003fast}Moenning, C. \& Dodgson, N. Fast marching farthest point sampling. (University of Cambridge, Computer Laboratory,2003)

\bibitem{rubner2000earth}Rubner, Y., Tomasi, C. \& Guibas, L. The earth mover's distance as a metric for image retrieval. {\em International Journal Of Computer Vision}. \textbf{40} pp. 99-121 (2000)
\bibitem{greenwade93}Greenwade, G. The Comprehensive Tex Archive Network (CTAN). {\em TUGBoat}. \textbf{14}, 342-351 (1993)
\bibitem{gao2015hypoelliptic}Gao, T. Hypoelliptic diffusion maps and their applications in automated geometric morphometrics. (Duke University,2015)
\bibitem{singercoifman}Singer, A. \& Coifman, R. Non-linear independent component analysis with diffusion maps. {\em Applied And Computational Harmonic Analysis}. \textbf{25} (2008)
\bibitem{shan2019probabilistic}Shan, S. Probabilistic Models on Fibre Bundles. (Duke University,2019)
\bibitem{steenrod1951topology}Steenrod, N. The topology of fiber bundles, volume 14 of. {\em Princeton Mathematical Series}. \textbf{16} (1951)
\bibitem{ravier2018eyes}Ravier, R. Eyes on the Prize: Improved Biological Surface Registration via Forward Propagation. {\em ArXiv Preprint ArXiv:1812.10592}. (2018)
\bibitem{puente2013distances}Puente, J. \& Others Distances and algorithms to compare sets of shapes for automated biological morphometrics. (Princeton, NJ: Princeton University,2013)
\bibitem{gao2021diffusion}Gao, T. The diffusion geometry of fibre bundles: Horizontal diffusion maps. {\em Applied And Computational Harmonic Analysis}. \textbf{50} pp. 147-215 (2021)
\bibitem{shan2019ariadne}Shan, S., Kovalsky, S., Winchester, J., Boyer, D. \& Daubechies, I. ariaDNE: A robustly implemented algorithm for Dirichlet energy of the normal. {\em Methods In Ecology And Evolution}. \textbf{10}, 541-552 (2019)
\bibitem{ravier2018algorithms}Ravier, R. Algorithms with Applications to Anthropology. (Duke University,2018)
\bibitem{he2004locality}He, X. \& Niyogi, P. Locality preserving projections. {\em Advances In Neural Information Processing Systems}. pp. 153-160 (2004)
\bibitem{roweis2000nonlinear}Roweis, S. \& Saul, L. Nonlinear dimensionality reduction by locally linear embedding. {\em Science}. \textbf{290}, 2323-2326 (2000)
\bibitem{sober2020manifold}Sober, B. \& Levin, D. Manifold approximation by moving least-squares projection (MMLS). {\em Constructive Approximation}. \textbf{52}, 433-478 (2020)
\bibitem{coifman2005geometric}Coifman, R., Lafon, S., Lee, A., Maggioni, M., Nadler, B., Warner, F. \& Zucker, S. Geometric diffusions as a tool for harmonic analysis and structure definition of data: Diffusion maps. {\em Proceedings Of The National Academy Of Sciences}. \textbf{102}, 7426-7431 (2005)
\bibitem{pearson1901l}Pearson, K. On lines and planes of closest fit to systems of points in space. {\em The London, Edinburgh, And Dublin Philosophical Magazine And Journal Of Science}. \textbf{2}, 559-572 (1901)
\bibitem{cox2000multidimensional}Cox, T. \& Cox, M. Multidimensional Scaling. (Chapman,2000)
\bibitem{faigenbaum2020manifold}Faigenbaum-Golovin, S. \& Levin, D. Manifold reconstruction and denoising from scattered data in high dimension. {\em Journal Of Computational And Applied Mathematics}. \textbf{421} pp. 114818 (2023)
\bibitem{tenenbaum2000global}Tenenbaum, J., De Silva, V. \& Langford, J. A global geometric framework for nonlinear dimensionality reduction. {\em Science}. \textbf{290}, 2319-2323 (2000)
\bibitem{shan2022diffusion}Shan, S. \& Daubechies, I. Diffusion maps: Using the semigroup property for parameter tuning. {\em ArXiv Preprint ArXiv:2203.02867}. (2022)
\bibitem{belkin2008towards}Belkin, M. \& Niyogi, P. Towards a theoretical foundation for Laplacian-based manifold methods. {\em Journal Of Computer And System Sciences}. \textbf{74}, 1289-1308 (2008)
\bibitem{belkin2006convergence}Belkin, M. \& Niyogi, P. Convergence of Laplacian eigenmaps. {\em Advances In Neural Information Processing Systems}. \textbf{19} (2006)
\bibitem{fisher1936use}Fisher, R. The use of multiple measurements in taxonomic problems. {\em Annals Of Eugenics}. \textbf{7}, 179-188 (1936)
\bibitem{boyer2011algorithms}Boyer, D., Lipman, Y., St. Clair, E., Puente, J., Patel, B., Funkhouser, T., Jernvall, J. \& Daubechies, I. Algorithms to automatically quantify the geometric similarity of anatomical surfaces. {\em Proceedings Of The National Academy Of Sciences}. \textbf{108}, 18221-18226 (2011)
\bibitem{daubechies2011continuous}Daubechies, Y. \& Others The continuous Procrustes distance between two surfaces. {\em ArXiv Preprint ArXiv:1106.4588}. (2011)

\bibitem{coifman2006diffusion}Coifman, R. \& Lafon, S. Diffusion maps. {\em Applied And Computational Harmonic Analysis}. \textbf{21}, 5-30 (2006)

\bibitem{gao2015hypoelliptic}Gao, T. Hypoelliptic diffusion maps and their applications in automated geometric morphometrics. (Duke University,2015)


\bibitem{bunn2011comparing}Bunn, J., Boyer, D., Lipman, Y., St. Clair, E., Jernvall, J. \& Daubechies, I. Comparing Dirichlet normal surface energy of tooth crowns, a new technique of molar shape quantification for dietary inference, with previous methods in isolation and in combination. {\em American Journal Of Physical Anthropology}. \textbf{145}, 247-261 (2011)
\bibitem{winchester2016morphotester}Winchester, J. MorphoTester: an open source application for morphological topographic analysis. {\em PloS One}. \textbf{11}, e0147649 (2016)
\bibitem{pampush2016introducing}Pampush, J., Winchester, J., Morse, P., Vining, A., Boyer, D. \& Kay, R. Introducing molaR: a new R package for quantitative topographic analysis of teeth (and other topographic surfaces). {\em Journal Of Mammalian Evolution}. \textbf{23} pp. 397-412 (2016)

\bibitem{winchester2014dental}Winchester, J., Boyer, D., St. Clair, E., Gosselin-Ildari, A., Cooke, S. \& Ledogar, J. Dental topography of platyrrhines and prosimians: convergence and contrasts. {\em American Journal Of Physical Anthropology}. \textbf{153}, 29-44 (2014)

\bibitem{shan2019ariadne}Shan, S., Kovalsky, S., Winchester, J., Boyer, D. \& Daubechies, I. ariaDNE: A robustly implemented algorithm for Dirichlet energy of the normal. {\em Methods In Ecology And Evolution}. \textbf{10}, 541-552 (2019)

\bibitem{fulwood2021reconstructing}Fulwood, E., Shan, S., Winchester, J., Gao, T., Kirveslahti, H., Daubechies, I. \& Boyer, D. Reconstructing dietary ecology of extinct strepsirrhines (Primates, Mammalia) with new approaches for characterizing and analyzing tooth shape. {\em Paleobiology}. \textbf{47}, 612-631 (2021)


\bibitem{herrera2016phylogeny}Herrera, J. \& Dávalos, L. Phylogeny and divergence times of lemurs inferred with recent and ancient fossils in the tree. {\em Systematic Biology}. \textbf{65}, 772-791 (2016)

\bibitem{wang2021statistical}Wang, B., Sudijono, T., Kirveslahti, H., Gao, T., Boyer, D., Mukherjee, S. \& Crawford, L. A statistical pipeline for identifying physical features that differentiate classes of 3D shapes. {\em The Annals Of Applied Statistics}. \textbf{15}, 638-661 (2021)

\bibitem{gao2018development}Gao, T., Yapuncich, G., Daubechies, I., Mukherjee, S. \& Boyer, D. Development and assessment of fully automated and globally transitive geometric morphometric methods, with application to a biological comparative dataset with high interspecific variation. {\em The Anatomical Record}. \textbf{301}, 636-658 (2018)

\bibitem{boyer2012earliest}Boyer, D., Costeur, L. \& Lipman, Y. Earliest record of Platychoerops (primates, Plesiadapidae), a new species from Mouras quarry, Mont de Berru, France. {\em American Journal Of Physical Anthropology}. \textbf{149}, 329-346 (2012)


\bibitem{rolfe2021slicermorph}Rolfe, S., Pieper, S., Porto, A., Diamond, K., Winchester, J., Shan, S., Kirveslahti, H., Boyer, D., Summers, A. \& Maga, A. SlicerMorph: An open and extensible platform to retrieve, visualize and analyse 3D morphology. {\em Methods In Ecology And Evolution}. \textbf{12}, 1816-1825 (2021)



\bibitem{watanabe2018many}Watanabe, A. How many landmarks are enough to characterize shape and size variation?. {\em PloS One}. \textbf{13}, e0198341 (2018)

\bibitem{warmlander2019landmark}Wärmländer, S., Garvin, H., Guyomarc'h, P., Petaros, A. \& Sholts, S. Landmark typology in applied morphometrics studies: What's the point?. {\em The Anatomical Record}. \textbf{302}, 1144-1153 (2019)


\bibitem{deng2022survey}Deng, B., Yao, Y., Dyke, R. \& Zhang, J. A Survey of Non-Rigid 3D Registration. {\em Computer Graphics Forum}. \textbf{41}, 559-589 (2022)
\bibitem{van2011survey}Van Kaick, O., Zhang, H., Hamarneh, G. \& Cohen-Or, D. A survey on shape correspondence. {\em Computer Graphics Forum}. \textbf{30}, 1681-1707 (2011)

\bibitem{ravier2018eyes}Ravier, R. Eyes on the Prize: Improved Biological Surface Registration via Forward Propagation. {\em ArXiv Preprint ArXiv:1812.10592}. (2018)

\bibitem{intanagonwiwat2003directed}Intanagonwiwat, C., Govindan, R., Estrin, D., Heidemann, J. \& Silva, F. Directed diffusion for wireless sensor networking. {\em IEEE/ACM Transactions On Networking}. \textbf{11}, 2-16 (2003)


\bibitem{al2013continuous}Al-Aifari, R., Daubechies, I. \& Lipman, Y. Continuous procrustes distance between two surfaces. {\em Communications On Pure And Applied Mathematics}. \textbf{66}, 934-964 (2013)


\bibitem{Gaussian_Process_Joint_Landmarking_of_Surfaces}Shan, S., \& Daubechies, I. Gaussian Process Joint
Landmarking for Collections of Surfaces (forthcoming).


\end{thebibliography}
